\documentclass[12pt,twoside,a4paper]{article}
\usepackage{mathrsfs}
\usepackage{mathrsfs}
\usepackage{mathrsfs}
\usepackage{mathrsfs}
\usepackage{mathrsfs}
\usepackage{mathrsfs}
\usepackage{mathrsfs}
\usepackage{stmaryrd}
\usepackage{amsfonts}
\usepackage{amsmath,amssymb}
\usepackage{mathrsfs}
\usepackage{hyperref}
\usepackage{graphicx}
\usepackage{graphicx}

\newtheorem{thm}{Theorem}[section]
\newtheorem{cor}[thm]{Corollary}
\newtheorem{lem}[thm]{Lemma}
\newtheorem{prop}[thm]{Proposition}

\numberwithin{equation}{section}\allowdisplaybreaks
\def\leq{\leqslant}
\def\ge{\geqslant}
\def\leq{\leqslant}
\def\geq{\geqslant}

\def\Real{{\mathbb{R}}}

\begin{document}

\title{\large\bf Necessary and sufficient conditions for the fractional Gagliardo-Nirenberg inequalities and applications to Navier-Stokes and generalized boson equations }
\author{\small \bf Hichem Hajaiej,$^{\dag}$ \ Luc Molinet,$^\S$  \  Tohru Ozawa$,^\sharp$   \  Baoxiang  Wang$^{\ddag} $ \\
{\small \it  $^\dag$Dept of Mathematics, King Saud University, P.O. Box 2455, 11451 Riyadh, Saudi Arabia}\\
{\small \it $^\S$L.M.P.T., Universit\'e Fran\c cois Rabelais Tours,  Parc Grandmont, 37200 Tours,  France }\\
{\small \it $^\sharp$Department of Applied Physics, Waseda University, Tokyo 169-8555, Japan}\\
{\small \it $^\ddag$LMAM, School of Mathematical Sciences, Peking University, Beijing 100871, China}\\
}\maketitle

\begin{abstract}

\rm Necessary  and  sufficient conditions for the generalized Gagliardo-Nirenberg  inequalities are obtained. For $0<q<\infty$, $0<p,  p_0,p_1 \leq \infty$, $s,s_0,s_1\in \mathbb{R}$,  $\theta \in (0,1)$,
\begin{align}
\| u\|_{\dot{B}_{p, q}^{s}}\lesssim  \| u\|_{\dot{B}_{p_0,
\infty}^{s_0}}^{1-\theta} \| u\|_{\dot{B}_{p_1, \infty}^{s_1}}^{
\theta} \label{besovgn-3}
\end{align}
holds if and only if
$
n/p-s   = (1-\theta) (n/p_0-
s_0 )
+\theta (n/p_1 - s_1),  \
  s_0- n/p_0 \neq s_1- n/p_1,   \
  s\leq (1-\theta)s_0 + \theta s_1$, and $p_0=p_1$ if $s = (1-\theta)s_0 + \theta s_1$.  Applying this inequality, we show that the solution of the Navier-Stokes equation at finite blowup time $T_m$ has a concentration phenomena in the critical space $L^3(\mathbb{R}^3).$   Moreover, we consider the minimization problem for the variational problem
$$
M_c =\inf \left\{ E(u): \ \|u_i\|^2_2=c_i >0, \ i=1,..., L \right\},
$$
where
$$
E(u)=\frac{1}{2}\|u\|^2_{\dot H^s}- \int_{\mathbb{R}^{2n}} G(u(x))V(x-y)G(u(y)) dxdy
$$
for $u=(u_1,..., u_L)\in (H^s)^L $ and show that $M_c$ admits a radial and radially decreasing minimizer under suitable assumptions on $s,$ $G$ and $V$.


\end{abstract}

\newpage

\section{Introduction}

The Gagliardo-Nirenberg (GN) inequality  is a fundamental tool in
the study of nonlinear partial differential equations, which was
discovered by Gagliardo \cite{Ga58}, Nirenberg \cite{Niren59} (see also \cite{Lad69}) in some special
cases. Throughout this paper, we denote by $L^p:= L^p(\mathbb{R}^n)$ the Lebesgue space, $\|\cdot\|_p: =\|\cdot\|_{L^p}$.  $C>1$ will denote positive
universal constants, which can be different at different places.
$a\lesssim b$ stands for $a\leq C b$ for some constant $C>1$, $a\sim
b$ means that $a\lesssim b$ and $b\lesssim a$. We write $a\wedge b
=\min(a,b)$, $a\vee b =\max(a,b)$. The classical integer version of the GN inequality can be stated as follows (see \cite{Fri} for instance):
\begin{thm} \label{GN}
 Let $1\leq p, p_0, p_1 \leq \infty$,  $\ell, m\in \mathbb{N}\cup
 \{0\}$, $\ell < m$, $\ell/m \leq \theta \leq 1$,
 and
\begin{align}
\frac{n}{p}- \ell = (1-\theta)\frac{n}{p_0}+\theta\left(\frac{n}{p_1}- m \right).
\label{gn-1}
\end{align}
Then we have for all $u\in C^\infty_0(\mathbb{R}^n)$,
\begin{align}
\sum_{|\alpha|=\ell}\|\partial^\alpha u\|_{ p}\lesssim  \|
u\|_{ p_0}^{1-\theta} \sum_{|\alpha|=m} \| \partial^\alpha
u\|_{ p_1}^{ \theta},
\end{align}
where we further assume $\ell /m \leq \theta <1$ if $m-\ell -n/p_1$
is an integer.
\end{thm}

The classical proof of the GN inequality is based
on the global derivative analysis in $L^p$ spaces, whose proof is
rather complicated, cf. \cite{Fri, Guo}. On the basis of the harmonic analysis techniques, there are some recent works devoted to  generalizations of the GN inequality, cf. \cite{ BGX00, BBM1, BBM2,CDPX01, CDDV, CMO98, EsVe97, Fri, Guo, GW, Led03, MaOz02, MeRi, Oru98, Strz06}.

In the first part of this paper, we consider the GN inequality with fractional order derivatives. First, we introduce some function spaces which will be frequently used,  cf. \cite{Tr}.  We denote by $\dot H^s_p:= (-\triangle)^{s/2} L^p$ the Riesz potential space, $\dot H^s= \dot H^s_2$, $H^s= L^2\cap \dot H^s$ for any $s\ge 0$.
Let $\psi$ be a smooth cut-off function supported in the ball $\{\xi:  |\xi| \leq 2\}$,  $\varphi= \psi(\cdot)- \psi(2\,\cdot)$.
We write
$ \varphi_k (\xi)=\varphi (2^{-k} \xi), \  k\in \mathbb{Z}.$
We see that
\begin{align}\label{eq1.5.2}
\sum_{k\in \mathbb{Z}} \varphi_k(\xi)=1,  \quad \xi\in \mathbb{R}^n \setminus \{0\}.
\end{align}
We introduce the homogeneous dyadic decomposition operators
$\triangle_k=\mathscr{F}^{-1}\varphi_k \mathscr{F}, \quad  k\in \mathbb{Z}.$
Let
$-\infty<s<\infty, \;  1\leq p, q\leq \infty.$
The space $\dot B^s_{p,q}$ equipped with  norm
\begin{align}\label{eq1.5.6}
 \|f\|_{\dot
B^s_{p,q}}:=\bigg(\sum^\infty_{k=-\infty} 2^{k s q}\|\triangle_k f
\|^q_p \bigg)^{1/q}
\end{align}
 is said to be a homogeneous Besov space (a tempered distribution $f\in \dot B^s_{p,q}$ modulo polynomials). Let
\begin{align}\label{eq1.5.8}
-\infty <s<\infty, \quad 1\leq p< \infty, \quad 1\leq q\leq  \infty.
\end{align}
The space $\dot F^s_{p,q}$ equipped with  norm
\begin{align}
 \label{eq1.5.10}
\|f\|_{\dot F^s_{p,q}}:=\bigg\|\bigg(\sum^\infty_{k=-\infty} 2^{k
s q} |\triangle_k f |^q\bigg)^{1/q} \bigg \|_p
\end{align}
is said to be a homogeneous Triebel-Lizorkin space (a tempered distribution  $f\in \dot F^s_{p,q}$ modulo polynomials).

In this paper we will obtain necessary  and sufficient conditions for the GN inequality in homogeneous Besov spaces $\dot B^s_{p,q}$ and Triebel-Lizorkin spaces $\dot F^s_{p,q}$.
As a corollary, we obtain that the GN inequality
also holds in fractional Sobolev spaces $\dot H^s_p$. The fractional GN inequalities in Theorems \ref{Besov-GN}, \ref{Besov-GN1} and \ref{Trieb-GN} below cover all of the available GN inequalities in \cite{ BGX00, BBM1, BBM2,CDPX01, CDDV, CMO98, EsVe97, Fri, Guo, GW, Led03, MaOz02, MeRi, Oru98, Strz06}  for both integer and fractional versions. Moreover, our results below clarify how the third indices $q$ in $\dot B^s_{p,q}$ and $\dot F^s_{p,q}$ contribute the validity of the GN inequalities.  We have
\begin{thm} \label{Besov-GN}
Let $0<p, p_0, p_1, q, q_0, q_1   \leq \infty$, $s,s_0,s_1\in
 \mathbb{R}$,  $0  \leq \theta \leq 1$.
Then the fractional GN inequality of the following type
\begin{align}
\| u\|_{\dot{B}_{p, q}^{s}}\lesssim  \| u\|_{\dot{B}_{p_0,
q_0}^{s_0}}^{1-\theta} \| u\|_{\dot{B}_{p_1, q_1}^{s_1}}^{\theta}
\label{besovgn-1}
\end{align}
holds for all $u\in \dot{B}_{p_0,
q_0}^{s_0} \cap \dot{B}_{p_1, q_1}^{s_1}$ if and only if
\begin{align}
\frac{n}{p}- s  & = (1-\theta)\left(\frac{n}{p_0}-
 s_0 \right)
+\theta\left(\frac{n}{p_1}- s_1 \right), \label{besob-200}\\
      s & \leq (1-\theta) s_0+ \theta s_1, \label{besovgn-2}\\
      \frac{1}{q}  & \leq   \frac{1-\theta }{q_0} + \frac{ \theta}{q_1}, \ \mbox{ if } \  p_0\neq p_1 \ and \  s= (1-\theta) s_0+ \theta s_1,    \label{besovgn-201}\\
      s_0  & \neq    s_1  \ \mbox{ or } \  \frac{1}{q}    \leq   \frac{1-\theta }{q_0} + \frac{ \theta}{q_1}, \ \mbox{ if } \  p_0= p_1 \ and \  s= (1-\theta) s_0+ \theta s_1,     \label{besovgn-202}\\
s_0-\frac{n}{p_0}  & \neq    s- \frac{n}{p}  \ \mbox{ or } \  \frac{1}{q}    \leq   \frac{1-\theta }{q_0} + \frac{ \theta}{q_1}, \ \mbox{ if }   \  s< (1-\theta) s_0+ \theta s_1.     \label{besovgn-203}
\end{align}
\end{thm}

\begin{thm} \label{Besov-GN1}
Let $0<q <\infty$, $0< p, \; p_0, \; p_1\leq \infty$, $0<\theta<1$,
$ s, s_0,  s_1\in \mathbb{R}$. Then the fractional GN inequality of the following type
\begin{align}
\| u\|_{\dot{B}_{p, q}^{s}}\lesssim  \| u\|_{\dot{B}_{p_0,
\infty}^{s_0}}^{1-\theta} \| u\|_{\dot{B}_{p_1, \infty}^{s_1}}^{
\theta} \label{besovgn-3}
\end{align}
holds if and only if
\begin{align}
\frac{n}{p}- s  & = (1-\theta)\left(\frac{n}{p_0}-
 s_0 \right)
+\theta\left(\frac{n}{p_1}- s_1 \right), \label{bespl0}\\
 s_0  -\frac{n}{p_0}& \neq s_1- \frac{n}{p_1},   \label{bespl}\\
 s & \leq  (1-\theta)s_0 + \theta s_1, \label{besovgn-4}\\
       p_0 & = p_1     \ \mbox{ if }   \  s= (1-\theta) s_0+ \theta s_1.     \label{besovgn-2022}
\end{align}
\end{thm}

In homogeneous  Triebel-Lizorkin spaces $\dot F^s_{p,q}$, we have the following

\begin{thm} \label{Trieb-GN}
Let $0<p, p_i, q  < \infty$, $s,s_0,s_1\in
 \mathbb{R}$,
  $0   < \theta < 1$. Then the fractional GN inequality of the following type
\begin{align}
\| u\|_{\dot{F}_{p, q}^{s}}\lesssim  \| u\|_{\dot{F}_{p_0,
\infty}^{s_0}}^{1-\theta} \| u\|_{\dot{F}_{p_1,
\infty}^{s_1}}^{\theta} \label{trieb-1}
\end{align}
holds if and only if
\begin{align}
\frac{n}{p}- s  & = (1-\theta)\left(\frac{n}{p_0}-
 s_0 \right)
+\theta\left(\frac{n}{p_1}- s_1 \right), \label{trieb-2}\\
& s\leq (1-\theta) s_0 + \theta s_1,  \label{trieb-3}\\
&  s_0\neq s_1 \mbox{ \ if \ } s= (1-\theta) s_0 + \theta s_1.  \label{trieb-4}
\end{align}
\end{thm}

 The following is the GN inequality with fractional derivatives.
\begin{cor} \label{Riesz-GN}
Let $1<p, p_0, p_1< \infty$, $s,\ s_1\in
 \mathbb{R}$,  $0  \leq \theta \leq 1$. Then the fractional GN inequality of the following type
\begin{align}
\| u\|_{\dot{H}_{p}^{s}}\lesssim  \| u\|_{L^{p_0}}^{1-\theta} \|
u\|_{\dot{H}_{p_1}^{s_1}}^{\theta} \label{riesz-1}
\end{align}
holds if and only if
\begin{align}
\frac{n}{p}- s = (1-\theta) \frac{n}{p_0}
+\theta\left(\frac{n}{p_1}- s_1 \right), \quad s \leq \theta
s_1. \label{riesz-2}
\end{align}
\end{cor}

\medskip

We will prove Theorems \ref{Besov-GN}--\ref{Trieb-GN} in Section 2. Relations with available GN inequalities are discussed in Section 3.  We remark that analogous results to Theorems \ref{Besov-GN}--\ref{Trieb-GN} and Corollary \ref{Riesz-GN}  also hold if one replaces all of the homogeneous spaces $\dot B^s_{p,q}, \, \dot F^s_{p,q}, \, \dot H^s_p$ by corresponding non-homogeneous spaces $ B^s_{p,q}, \, F^s_{p,q}, \, H^s_p$, respectively. We will list those results in Section \ref{nonhfs}.

In the second part of this paper we consider some applications of the fractional GN inequality. First, We study the Cauchy problem for the Navier-Stokes (NS) equation
\begin{align}\label{eq4.0.1}
 u_t  - \Delta u + (u\cdot \nabla) u + \nabla p=0, \ \ {\rm div} \, u=0,  \quad u(0,x)=u_0(x), \index{NS}
\end{align}
where $\Delta= \sum^n_{i=1}\partial^2_{x_i}$,
$\nabla=(\partial_{x_1},..., \partial_{x_n})$, ${\rm div}\, u
=\partial_{x_1}u_1+...+ \partial_{x_n}u_n$,  $u=(u_1,...,u_n)$ and
$p$ are real-valued unknown functions of $(t,x) \in [0,T_m) \times \mathbb{R}^n$ for some $T_m>0$,
$u_0=(u^1_0,...,u^n_0)$ denotes the initial value of $u$ at $t=0$. It is known that NS equation is local well posed in $L^n$, namely, for initial data $u_0\in L^n(\mathbb{R}^n)$, there exists a unique local solution $u\in C([0,T_m); L^n) \cap L^{2+n}_{\rm loc}(0,T_m; L^{2+n})$ (cf. \cite{Kato84, KeKo09}). Whether the local solution can be extended to a global one is still open.   Recently,  Escauriaza,  Seregin and  ${\rm \check{S}}$ver\'{a}k \cite{EsSeSv03}
showed that any ``Leray-Hopf" weak solution in 3D which remains bounded in $L^3(\mathbb{R}^3)$
cannot develop a singularity in finite time. Kenig and Koch \cite{KeKo09} gave an alternative approach to this problem by substituting $L^3$ with $\dot H^{1/2}$.  Dong and Du \cite{DoDu09} generalized their results in higher spatial dimensions $n\ge 3$. Noticing that $L^3 \subset B^{-1}_{\infty, \infty}$ in 3D is a sharp embedding,  for any solution $u$ of the NS equation  in $C([0, T^*); L^3)$,   we see that $u\in C([0, T^*);  B^{-1}_{\infty, \infty}) $.  May \cite{May03} (see also \cite{LeRi02})
prove that if $T^*<\infty$, then there exists a constant $c>0$ independent of the solution of NS equation such that
$
\lim\sup_{t\to T^*} \|u(t)-\omega \|_{B^{-1}_{\infty, \infty}} \ge c
$
for all $\omega \in \mathscr{S}$. In this paper we will use the fractional GN inequality to study the finite time blowup solution and we have the following concentration result:

\begin{thm} \label{NS concent}
Let $n=3$ and $u \in C([0,T_m); L^n \cap L^2) \cap L^{2+n}_{\rm loc}(0,T_m; L^{2+n})$ be the solution of NS equation with maximal existing time $T_m < \infty$.  Then there exist $c_0>0$ and $\delta>0$ such that
\begin{align}
\overline{\lim_{t\nearrow T_m}} \sup_{x_0 \in \mathbb{R}^n}  \int_{|x-x_0| \le (T_m-t)^\delta} |u(t, \, x-x_0) |^n dx \ge c_0,  \label{concent}
\end{align}
where the constant $c_0>0$ only depends on $\|u_0\|_n$, $\delta$ can be chosen as any positive constant less than $2/n^2$.
\end{thm}

As the second application of fractional GN inequalities, we consider the existence of the radial and radially decreasing non-negative solutions for the following  system:
\begin{align}
(m^2-\triangle)^s u_i - \left[ G(u) * V \right] \partial_{i} G(u) + r_i u_i =0, \  \ i=1,...,L ,   \label{FLE}
\end{align}
where $m^2\geq 0$, $u=(u_1,..., u_L) $,  $u_i \ge 0$ and $u \neq 0$,  $G: \ \mathbb{R}^L_+ \to \mathbb{R}_+=[0, \infty) $ is a differentiable function, $\partial_i G (v_1,..., v_L) :=  \partial G(v_1,..., v_L) /\partial v_i  $. $V (x) = |x|^{-(n-\beta)}$, $*$ denotes the convolution in $\mathbb{R}^n$, $r_i>0$. In order to work out a desired solution of \eqref{FLE}, it suffices to consider the existence of the radial and radially decreasing non-negative and non-zero minimizers of the following variational problem. We write for $c_1,..., c_L >0$,
\begin{align}
S_c  =  \left\{u= (u_1,..., u_L) \in (H^s)^L : \  \|u_i\|^2_2 =c_i, \ i=1,...,L \right\}.  \label{mass}
\end{align}
We will consider the variation problem
\begin{align}
M_c = \inf \{E(u): \ u \in S_c, \ c_1,..., c_L>0 \},  \label{ques-minizer}
\end{align}
where
\begin{align}
E(u) = \frac{1}{2} \sum^L_{i=1} \|(m^2+|\xi|^2)^{s/2}\widehat{u}_i\|^2_{2} -\int \int G( u  (x) ) V(|x-y|) G( u(y) ) dxdy. \label{energy}
\end{align}

Fractional calculus has gained tremendous popularity during the last two decades thanks to its applications in widespread domains of sciences, economics and engineering, see \cite{AgBeNiOu, BaDiScTr, KiSrTr, La}. Fractional powers of the Laplacian arise in many areas. Some of the fields of applications of fractional Laplacian models include  medicine where the equation of motion of semilunar heart value vibrations and stimuli of neural systems are modeled by a Capulo fractional Laplacian; cf. \cite{Sh, LuHiSpFa}. It also appears in modeling populations \cite{RiTrVaVa}, flood  flow, material viscoelastic theory, biology dynamics, earthquakes, chemical physics, electromagnetic theory, optic, signal processing, astrophysics, water wave, bio-sciences dynamical process and turbulence; cf. \cite{AgBeNiOu, AbBoFeSa, BaDiScTr, BoGe, CaVa, CaCoGaOr, Co, FrSe, FrLe, KiSrTr, La, LiYau,  MaMcTa, Ma, TaZa}.

In \cite{LiYau}, Lieb and Yau studied the existence and symmetry of ground state solutions for the boson equation in three dimensions:
\begin{align}
(m^2-\triangle)^{1/2} u  - \left( |x|^{-1} * u^2  \right) u +r u=0,     \label{Hartree1}
\end{align}
Taking $G(u)=u^2$ and $V(x)=|x|^{-1}$ in three dimensions, \eqref{FLE} is reduced to \eqref{Hartree1}.  The variational problem associated with \eqref{Hartree1} is
\begin{align}
M^{(3)}_c = \inf \left\{ \frac{1}{2}\|(m^2+|\xi|^2)^{1/4}\widehat{u}\|^2_{2}- \int_{\mathbb{R}^3}\int_{\mathbb{R}^3}  \frac{|u (x)|^2  |u (y)|^2 }{|x-y|} dxdy: \ u \in H^{1/2}, \  \|u\|^2_2=c\right\}. \label{Hartree2}
\end{align}
As indicated in \cite{LiYau}, \eqref{Hartree1} and \eqref{Hartree2}  play  a fundamental role in the mathematical theory of gravitational collapse of boson stars. Indeed, Lieb and Yau essentially showed that for $s=1/2$, there exists $c_*>0$, such that  \eqref{Hartree1} has a non-negative radial solution  if and only if $c= c_*$. It was proven in \cite{LiYau} that boson stars with total
mass strictly less than $c^*$ are gravitationally stable, whereas boson stars whose total mass exceed $c^*$
  may undergo a ``gravitational collapse"  based on variational arguments and many-body quantum
theory. The main tools used by Lieb and Yau are the Hardy-Littlewood-Sobolev inequality together with some rearrangement inequalities. Inspired and motivated by Lieb and Yau's work, Frank and Lenzemann \cite{FroLe} recently showed the uniqueness of ground states to \eqref{Hartree1} in 1D.

Taking $G(u)=u^2$ and $V(x)= |x|^{-(n-2)}$ in $n$-dimensions with $n\ge 3$, \eqref{FLE} is reduced to the general Choquard-Peckard equation
\begin{align}
(m^2-\triangle)^{s} u  - \left( |x|^{- (n-2)} * u^2  \right) u +r u =0.     \label{ChPe1}
\end{align}
The variational problem associated with \eqref{ChPe1} is
\begin{align}
M^{(n)}_c = \inf \left\{ \frac{1}{2}\|(m^2+|\xi|^2)^{s/2}\widehat{u}\|^2_{2}- \Upsilon_2 (u): \, u\in H^s, \ \|u\|^2_2=c\right\},   \label{ChPe2}
\end{align}
where
\begin{align}
\label{Upsilon}
\Upsilon_\beta (u) =  \int   \frac{|u (x)|^2  |u (y)|^2 }{|x-y|^{n-\beta} }dxdy.
\end{align}
Taking $G(u)=u^2_1+ u^2_2$ and $V(x)= |x|^{-1}$ in $3$-dimensions, \eqref{FLE} is reduced to the following system
\begin{align}
(m^2-\triangle)^{s} u_i - \left( |x|^{- 1} * (u^2_1+ u^2_2)  \right) u_i +r_i u_i =0,  \ i=1,2,   \label{system1}
\end{align}
which was studied in \cite{AsSq} and \cite{FroLe} in the cases $s=1$ and $s=1/2$, respectively.  If we treat $u=(u_1,u_2)$ and $\|u\|^2_X = \|u_1\|^2_X+ \|u_2\|^2_X $, we see that the   variational problem associated with \eqref{system1}  is the same as in \eqref{ChPe2} if one constraint $\|u_1\|^2_2 + \|u_2\|^2_2=c$ is considered.

\medskip

Now we state our main result on the existence of the minimizer of \eqref{ques-minizer}. There are two kinds of basic nonlinearities, one is $G(u)= u_1^{\mu_1}  ...  u^{\mu_L}_L$ and another is  $G(u)= u_1^\mu +...+ u^\mu_L$. For the former case, we need to use $m$-constraints $\|u_i\|^2_2= c_i>0$ to prevent the situation that the second term of $E(u)$ in \eqref{energy} vanishes. For the later case, one can use  one constraint $\|u_1\|^2_2+...+ \|u_L\|^2_2=c$. Let $ s\ge (n-\beta)/2.$ We first consider the former case and our main assumptions on $G$ are the following:
\begin{itemize}
\item[\rm (G1)]  $G: \mathbb{R}^L_+ \ni (v_1,..., v_L)   \to G(v_1,..., v_L) \in \mathbb{R}_+$ is a continuous function and there  exists  $ \mu \in [2,  1 + (2s+\beta)/n)$ such that
\begin{align}
  G(v) \leq C(|v|^2+ |v|^\mu), \ v=(v_1,..., v_L) .   \label{growthG}
\end{align}
Moreover, there exist  $\alpha_i >0$ such that for all $0<v_1,...,v_L \ll 1$,
\begin{align}
  G(v) \ge  c v^{\alpha_1}_1 v^{\alpha_2}_2... v^{\alpha_L}_L.   \label{growthG2}
\end{align}
where $ 0 < n+\beta-n(\alpha_1+...+\alpha_L) +2s$.

\item[\rm (G2)] If $v$ has a zero component, then  $G(v)=0$. The function $G\otimes G:  \mathbb{R}^L_+ \times \mathbb{R}^L_+ \ni (u,v) \to  G(u)G(v) \in \mathbb{R}_+$ is a super-modular\footnote{$F$ is said to be a supermodular if (\cite{Lor})
$$F(y + h e_i + ke_j ) + F(y) \geq  F(y + he_i ) + F(y + k e_j)  \  (i\neq j, \ h, k > 0), $$
where $y = (y_1, . . . , y_L)$, and $e_i$ denotes the {\it i}-th standard basis vector in $\mathbb{R}^L$. It is known that a smooth
function is a supermodular  if all its mixed second partial derivatives are nonnegative.}.

\item[\rm (G3)]  $G(t_1v_1,...,t_Lv_L) \ge t_{\max } G(v_1,..., v_L) $ for any $t_i  \ge 1$, where $t_{\max} = \max (t_1,...,t_L)$.

\end{itemize}
Noticing that $v^{\alpha_1}_1 v^{\alpha_2}_2... v^{\alpha_L}_L \le |v|^{\alpha_1+...+\alpha_L}$, we see that condition \eqref{growthG} covers the nonlinearity $G(v)=v^{\alpha_1}_1 v^{ \alpha_2}_2... v^{\alpha_L}_L$ if $\alpha_1+...+\alpha_L \in [2,\mu]$.  Our main result on the existence of the minimizer of \eqref{ques-minizer}  is the following:

\begin{thm}\label{exist-minimizer}
Let $m^2\geq 0$, $0<\beta <n$, $s>(n-\beta)/2 $.  Assume that conditions {\rm  (G1)--(G3)} are satisfied. Then \eqref{ques-minizer} admits a radial and radially decreasing minimizer in $(H^s)^L $.
\end{thm}
We point out that both conditions $s\geq (n-\beta)/2$ and $ 0 < n+\beta-n(\alpha_1+...+\alpha_L) +2s$ are necessary for Theorem \ref{exist-minimizer}. Indeed, we can give a counterexample to show that $M_c=-\infty$ if $s<(n-\beta)/2$ or $ 0 > n+\beta-n(\alpha_1+...+\alpha_L) +2s$ for a class of nonlinearities $G(u)$.

The endpoint case $s= (n-\beta)/2$ can not be handled in Theorem \ref{exist-minimizer}. Note that for $s= (n-\beta)/2$,  we have $\mu=2$ in \eqref{growthG}, a basic example is $G(u)= u_1^2+...+ u^2_L$. Now we consider  the variational problem
 \begin{align}
M^{(n)}_{c, \beta} = \inf \left\{ \frac{1}{2}\|(m^2+|\xi|^2)^{s/2}\widehat{u}\|^2_{2}- \Upsilon_\beta (u): \, u\in (H^s)^L, \  \|u\|^2_2=c >0\right\}. \label{ChPe200}
\end{align}
where $u=(u_1,..., u_L) $, $|u|^2=u^2_1+...+u^2_L$  and $\|u\|^2_X = \|u_1\|^2_X+...+ \|u_L\|^2_X$.      Using the definition of the Riesz potential, the Plancherel  identity,  the Hardy-Littlewood-Sobolev,  and  fractional GN inequalities, we have
\begin{align}\label{Upsilon}
\Upsilon_\beta (u) & = C(n, \beta) \int |u (x)|^2 [(-\Delta)^{-\beta/2} |u|^2] (x) dx =\|(-\Delta)^{-\beta/4} |u|^2\|^2_2  \nonumber\\
 &  \leq C \left(\|u_1\|^2_{ 4n/(n+\beta)} +...+ \|u_L\|^2_{4n/(n+\beta)} \right)^2 \nonumber\\
 &\leq C \left(\|u_1\|_{2} \|u_1\|_{\dot H^{(n-\beta)/2}}  +...+ \|u_L\|_{2} \|u_L\|_{\dot H^{(n-\beta)/2}}  \right)^2 \nonumber\\
  & \leq C  \|u\|^2_2 \|u\|^2_{\dot H^{(n-\beta)/2}}.
\end{align}
  Define
\begin{align}\label{C*}
 C^* =\sup_{u\in H^{(n-\beta)/2} \setminus \{ 0\} } \frac{\Upsilon_\beta (u)}{\|u\|^2_2 \|u\|^2_{\dot H^{(n-\beta)/2}}}.
\end{align}

\begin{thm}\label{exist-minimizer2}
Let $m^2= 0$, $0<\beta <n$,  $s= (n-\beta)/2$,  $G(u)=u^2_1+...+u^2_L$.   Then \eqref{ChPe200} admits a radial and radially decreasing minimizer in $ (H^s)^L  $ if and only if $c=1/2C^*$.
\end{thm}
As a straightforward consequence of Theorem \ref{exist-minimizer}, we see that \eqref{ChPe2}
admits a radial and radially decreasing minimizer in $ H^{(n-2)/2}$ if and only if $c=1/2C^*$, where $\beta=2$ in the definition of $C^*$.

\medskip

In the case $m^2>0$ we have the following
\begin{thm}\label{exist-minimizer4}
Let  $m^2>0$, $0<\beta<n$, $s= (n-\beta)/2$, $c>0$.   Then we have
\begin{itemize}
\item[\rm (1)] If $n>2+\beta$, then \eqref{ChPe200} has no minimizer in in $ (H^s)^L  $.

\item[\rm (2)] If $n<2+\beta$, then \eqref{ChPe200} admits a radial and radially decreasing minimizer in $ (H^s)^L  $ if and only if $0<c<1/2C^*$.

 \item[\rm (3)] If $n=2+\beta$, then \eqref{ChPe200} admits a radial and radially decreasing minimizer in $ (H^s)^L  $ if and only if $ c= 1/2C^*$.

 \end{itemize}
\end{thm}

\section{Proofs of the GN inequalities}

The following is an interpolation inequality in Besov spaces, which is very useful in nonlinear estimates, see \cite{GV1,GW}.

\begin{prop}[Convexity H\"{o}lder's inequality]\label{convholder}
Let $0< p_i, q_i \leq \infty$,
$0\leq \theta_i\leq 1$, $\sigma_i, \sigma\in \mathbb{R}$ $ (i=1,\ldots,
N)$, $\sum^N_{i=1}\theta_i=1$,
$\sigma=\sum^N_{i=1}\theta_i\sigma_i,$
$1/p=\sum^N_{i=1}\theta_i/p_i,$   $1/q = \sum^N_{i=1}\theta_i/q_i$.
Then $\cap^N_{i=1}\dot B^{\sigma_i} _{p_i, q_i}\subset \dot
B^\sigma_{p,q}$  and for any  $v\in\cap^N_{i=1}\dot B^{\sigma_i} _{p_i,
q_i}$,
\begin{align*}
\|v\|_{\dot B^\sigma_{p,q}}\leq \prod^N_{i=1}\|v\|^{\theta_i}
_{\dot B^{\sigma_i}_{p_i, q_i}}.
\end{align*}
This estimate also holds if one substitutes $\dot B^\sigma_{p,q}$ by $\dot F^\sigma_{p,q}$ ($p,p_i\not= \infty$).
\end{prop}
In the convexity H\"older inequality, condition  $1/q = \sum^N_{i=1}\theta_i/q_i$ can be replaced by  $1/q \leq \sum^N_{i=1}\theta_i/q_i$. Indeed, noticing that $\ell^q \subset \ell^p$ for all $q\leq p$, we see that Proposition \ref{convholder} still holds if
 $1/q < \sum^N_{i=1}\theta_i/q_i$. In \cite{GV1,GW}, Proposition \ref{convholder}  was stated as the case $1\le p_i, q_i\le \infty$, however, the proof in \cite{GW} is also adapted to the case $0<p_i, q_i \le \infty$.

\medskip

\noindent {\bf Proof of Theorem \ref{Besov-GN}} (Sufficiency)  First, we consider the case $1/q   \leq (1-\theta )/ q_0  +   \theta/ q_1$.  By
\eqref{besovgn-2}, we have
\begin{align}
\frac{1}{p}- \frac{1-\theta}{p_0}-  \frac{\theta}{p_1} = \frac{s}{n}
- (1-\theta)\frac{s_0}{n}  - \theta \frac{s_1}{n}:= -\eta \leq 0.
\label{besovgn-5}
\end{align}
Take $p^*$ and $s^*$ satisfying
$$
\frac{1}{p^*}=\frac{1}{p} + \eta, \quad  s^*= s+ n \eta.
$$
Applying the convexity H\"older inequality, we have
\begin{align}
\|f\|_{\dot B^{s^*}_{p^*,q}}
 \leq  \|f\|_{\dot B^{s_0}_{p_0,q_0}}^{1-\theta} \|f\|_{\dot
B^{s_1}_{p_1,q_1}}^{\theta}.
 \label{besovgn-6}
\end{align}
Using the inclusion $\dot B^{s^*}_{p^*,q} \subset \dot B^{s}_{p,q}$,
we get the conclusion.

Next, we need to consider the following two cases: (i) $s=(1-\theta ) s_0+\theta s_1$, $p_0=p_1$ and $s_0\neq s_1$; (ii) $s< (1-\theta ) s_0+\theta s_1$ and $s-n/p\neq s_0- n/p_0$.   We can show that
\begin{align}
\|f\|_{\dot B^{s }_{p ,q}}
 \leq  \|f\|_{\dot B^{s_0}_{p_0,\infty}}^{1-\theta} \|f\|_{\dot
B^{s_1}_{p_1,\infty}}^{\theta},
 \label{besovgn-60}
\end{align}
see below, the proof of Theorem \ref{Besov-GN1}.  \eqref{besovgn-60} implies the result, as desired.

\medskip

 (Necessity) By scaling,
$$
\|f(\lambda\cdot)\|_{\dot B^{s}_{p,q}} \sim \lambda ^{s-n/p} \|f\|_{\dot
B^{s}_{p,q}}, \quad \lambda\in 2^{\mathbb{Z}}.
$$
Hence, if \eqref{besovgn-1} holds, then
$$
\lambda ^{s-n/p- [(1-\theta)(s_0-n/p_0)+\theta (s_1-n/p_1)]} \leq C.
$$
Letting $\lambda\to 0$ or $\lambda \to \infty$, we immediately obtain that $s-n/p-
[(1-\theta)(s_0-n/p_0)+\theta (s_1-n/p_1)=0$.

Next, we show that $s-s_0 \leq \theta (s_1-s_0)$. Assume on the contrary that $s-s_0
> \theta (s_1-s_0)$. Assume that $s_0=0$. Let $\varphi$ satisfy supp $\varphi \subset \{\xi: 1/2\leq
|\xi| \leq 3/2\}$ and $\varphi(\xi)=1$ for $3/4\leq |\xi| \leq 1$.
So,
 $\varphi(2^{-j} \xi)=1$ if $  3\cdot  2^{j-2}\leq |\xi| \leq 2^j.
 $
Denoting
\begin{align}
\rho_j (\xi) = \varphi (2(\xi- \xi^{(j)} )), \ \  \xi^{(j)} = (7\cdot 2^{j-3},0,...,0).
\label{conter}
\end{align}
and for sufficiently small $\varepsilon>0$, we write
\begin{align}
\hat{f} (\xi) = \sum^N_{j=100} 2^{\varepsilon j}\rho_j (\xi).
\end{align}
This leads to
$$
\|f\|^q_{\dot B^s_{p,q}} =\sum^N_{j=100} 2^{(s+\varepsilon)qj}
\|\mathscr{F}^{-1} (\varphi_j \rho_j)\|^q_p.
$$
Noticing that  $\varphi_j(\xi)=1$ for $\xi\in {\rm supp}\,\rho_j$, we have
$$
\|\mathscr{F}^{-1} (\varphi_j \rho_j)\|_p=\|\mathscr{F}^{-1}
\rho_j\|_p= \|\mathscr{F}^{-1} \rho_0\|_p.
$$
Hence,
$$
\|f\|_{\dot B^s_{p,q}} \sim 2^{(s+\varepsilon)N}.
$$
Similarly,
$$
\|f\|_{\dot B^0_{p_0,q_0}} \sim 2^{\varepsilon N}, \quad \|f\|_{\dot
B^{s_1}_{p_1,q_1}} \sim 2^{(s_1+\varepsilon) N}.
$$
By \eqref{besovgn-1}, we obtain that $2^{(s+\varepsilon)N} <
2^{\varepsilon N} 2^{s_1\theta N} $. However, for sufficiently large $N$, it contradicts  the fact
$s>\theta s_1$. Substituting $s$ by $s-s_0$, we get the proof in the case $s_0\not=0 $.

Thirdly, we consider the case $p_0\neq p_1$ and $s=(1-\theta)s_0 + \theta s_1$ and show that $1/q\leq (1-\theta)/q_0 + \theta/ q_1$.    Put
\begin{align}
\lambda = \frac{s_1-s_0}{n(1/p_0 -1/p_1)}. \label{lambda}
\end{align}
We see that
\begin{align} \label{homog}
s+ n\lambda \left( \frac{1}{p} -1\right) =  s_0 + n\lambda \left( \frac{1}{p_0} -1\right) =  s_1 + n\lambda \left( \frac{1}{p_1} -1\right).
 \end{align}
{\it Case 1.} We consider the case $\lambda \ge 0$. Let $\varphi$ and $\xi^{(j)}$ be as in \eqref{conter}.  Put
$$
\varrho^\lambda_j := \varphi (2^{\lambda j} (\xi- \xi^{(j)}))
$$
and
\begin{align} \label{F}
\widehat{F} =\sum^J_{j=100} 2^{-sj-n\lambda(1/p-1)j} \varrho^\lambda_j.
\end{align}
Since ${\rm supp} \ \widehat{F}$ overlaps only one ${\rm supp} \ \varphi_j$ for all $j\in \mathbb{Z}$ and for $j \ge 100$,
$$
\|\triangle_j \mathscr{F}^{-1} \varrho^\lambda_j\|_p = \|  \mathscr{F}^{-1} \varrho^\lambda_j\|_p  \sim  2^{n\lambda j(1/p-1)},
$$
we have
\begin{align} \label{normF}
 \|F\|^q_{\dot B^s_{p,q}} & =  \sum_{j\in \mathbb{Z}} (2^{s  j} \|\triangle_j F\|_p)^q \nonumber\\
  & =  \sum^J_{j=100}  (2^{s  j}  \|\triangle_j F\|_p)^q \nonumber\\
   & =  \sum^J_{j=100} ( 2^{-n\lambda  (1/p-1)j}  \|\triangle_j \mathscr{F}^{-1} \varrho^\lambda_j\|_p )^q \nonumber\\
  &    \sim J,
\end{align}
which means that $\|F\|_{\dot B^s_{p,q}} \sim J^{1/q}$.  On the other hand, in view of \eqref{homog} and \eqref{F}, we see that
\begin{align} \label{Fa}
\widehat{F} =\sum^J_{j=100} 2^{-s_0 j-n\lambda(1/p_0 -1)j} \varrho^\lambda_j =\sum^J_{j=100} 2^{-s_1 j-n\lambda(1/p_1 -1)j} \varrho^\lambda_j.
\end{align}
In an analogous way to \eqref{normF}, we find that
\begin{align} \label{normFa}
 \|F\|_{\dot B^{s_0}_{p_0,q_0}}
    \sim J^{1/q_0}, \ \ \|F\|_{\dot B^{s_1}_{p_1,q_1}}
    \sim J^{1/q_1},
\end{align}
By \eqref{besovgn-1}, we have $J^{1/q} \lesssim J^{(1-\theta)/q_0} J^{\theta/ q_1}$ for any $J\gg 1$.  It follows that $1/q \leq (1-\theta)/q_0 + \theta/ q_1$.

{\it Case 2}. We consider the case $\lambda <0$. Denote
$$
\varphi^{(N)} = \varphi(2^{-N} \, \cdot),  \ \  \varphi^{(N)}_j  = \varphi(2^{-j -N} \, \cdot), \ \ \triangle_{j,N} = \mathscr{F}^{-1}  \varphi^{(N)}_j \mathscr{F}.
$$
It is easy to see that
$$
\|f\|^{(N)}_{\dot B^s_{p,q}} = \left( \sum_j (2^{sj} \|\triangle_{j,N}\|_p)^q\right)^{1/q}
$$
is an equivalent norm on $\dot B^s_{p,q}$ (see also \cite{Tr}).  Let
\begin{align} \label{F-B}
\widehat{F} =\sum^J_{j=100} 2^{-sj-n\lambda(1/p-1)j}  \varphi (2^{\lambda j} \, \cdot).
\end{align}
Assuming that $N \ge 100  (|\lambda|+1)$, analogously to the above, we have from the definition of $\|\cdot\|^{(N)}_{\dot B^s_{p,q}} $ that
\begin{align} \label{normFa}
\|F\|^{(N)}_{\dot B^{s}_{p,q}}
    \sim J^{1/q}, \ \  \|F\|^{(N)}_{\dot B^{s_0}_{p_0,q_0}}
    \sim J^{1/q_0}, \ \ \|F\|^{(N)}_{\dot B^{s_1}_{p_1,q_1}}
    \sim J^{1/q_1}.
\end{align}
By \eqref{besovgn-1} we have $1/q \leq (1-\theta)/q_0 + \theta/ q_1$.

\medskip

Fourthly, we show the necessity of \eqref{besovgn-202}. If not, then we have $p_0=p_1=p$, $s_0=s_1=s$ and $1/q> (1-\theta)/q_0 + \theta /q_1$.  Let
\begin{align} \label{F-C}
\widehat{F} =\sum^J_{j=100} 2^{-sj+ n(1/p-1)j}  \varphi (2^{-j} \, \cdot).
\end{align}
We easily see that for $N\gg 1$,
\begin{align} \label{normFc}
\|F\|^{(N)}_{\dot B^{s}_{p,q}}    \sim J^{1/q}, \ \  \|F\|^{(N)}_{\dot B^{s }_{p ,q_0}}
    \sim J^{1/q_0}, \ \ \|F\|^{(N)}_{\dot B^{s }_{p ,q_1}}
    \sim J^{1/q_1}.
\end{align}
We have $1/q \leq (1-\theta)/q_0 + \theta/ q_1$, which is a contradiction.

\medskip

Finally, we show the necessity of \eqref{besovgn-203}. Assume for a contrary that $s-n/p =s_0-n/p_0$ and $1/q > (1-\theta)/q_0 + \theta /q_1.$  Using the same way as in \eqref{F-C} and \eqref{normFc}, we have a contraction.
$\hfill \Box$

\medskip
\medskip

\noindent {\bf Proof of Theorem \ref{Besov-GN1}.} (Sufficiency) We can assume that $s_0=0$ and
the case $s_0\not= 0$ can be shown by a similar way.

{\bf Step 1.}  We consider the case $p\ge p_0 \vee p_1$. By  definition,
\begin{align}
\| u\|_{\dot{B}_{p, q}^{s}}= \left(\sum\limits_{N\, {\rm
dyadic}}N^{s q}\|\triangle_N u\|_p^q \right)^{1/q}.
\end{align}
From \eqref{besovgn-4}, it follows that
\begin{align}
\theta\left(\frac{n}{p} -\frac{n}{p_1}+ s_1-s \right)
=(1-\theta)\left(  s +\frac{n}{p_0}- \frac{n}{p}\right).
\end{align}
Since $0<\theta<1$,   \eqref{bespl} implies that $ \big(\frac{n}{p}
-\frac{n}{p_1}+s_1-s\big) \big( s+\frac{n}{p_0}-
\frac{n}{p}\big)>0.$

{\it Case }1. We consider the case
\begin{align}
s_1-s +\frac{n}{p} -\frac{n}{p_1}>0, \quad  s+\frac{n}{p_0}-
\frac{n}{p}>0. \label{inter-6}
\end{align}
Using the inclusion $\dot B^s_{p, r_1} \subset \dot B^s_{p, r_2}$ for any $r_1\le r_2$, it suffices to consider the case $q<1/2, \; q^{-1}\in \mathbb{N}$. For brevity, we write
$K:=q^{-1}$.
\begin{align}
\| u\|_{\dot{B}_{p, q}^{s}} & \leq   \sum\limits_{N_1\geq ... \geq N_K} \left( N_1^{s }... N_K^{s
}\|\triangle_{N_1} u\|_p...
\|\triangle_{N_K} u\|_p \right)^{q^2} \nonumber\\
& \quad \times \left (N_1^{s }\ldots N_K^{s }\|\triangle_{N_1}
u\|_p\ldots \|\triangle_{N_K} u\|_p \right )^{q(1-q)}.
\label{inter-7}
\end{align}
In view of  Bernstein's inequality,
\begin{align}
 \|\triangle_{N} u\|_p\leq
N^{\frac{n}{p_0}-\frac{n}{p}}\|\triangle_{N} u\|_{p_0}, \ \  \|\triangle_{N} u\|_p\leq
N^{\frac{n}{p_1}-\frac{n}{p}}\|\triangle_{N} u\|_{p_1}.
\end{align}
We can choose $a \in (0,1], k \geq 1$ satisfying $\theta K = k-1+a$. Hence,
\begin{align}
\|\triangle_{N_1}&  u\|_p... \| \triangle_{N_K} u\|_p\nonumber\\
 =& ( \|\triangle_{N_1} u\|_p...
\|\triangle_{N_{k-1}} u\|_p\|\triangle_{N_{k}} u\|_p^a)
(\|\triangle_{N_{k}} u\|_p^{1-a}
\|\triangle_{N_{k+1}} u\|_p\ldots \|\triangle_{N_{K}} u\|_p) \nonumber\\
\lesssim &\; N_k^{(1-a)(\frac{n}{p_0}-\frac{n}{p})}
N_{k+1}^{\frac{n}{p_0}-\frac{n}{p}}\ldots
N_K^{\frac{n}{p_0}-\frac{n}{p}}\|\triangle_{N_k}
u\|_{p_0}^{1-a}\|\triangle_{N_{k+1}}
u\|_{p_0}\ldots\|\triangle_{N_K} u\|_{p_0}\nonumber\\
&  \times  N_1^{\frac{n}{p_1}-\frac{n}{p}}\ldots
N_{k-1}^{\frac{n}{p_1}-\frac{n}{p}}N_k^{a(\frac{n}{p_1}-\frac{n}{p})}
\|\triangle_{N_1} u\|_{p_1}\ldots \|\triangle_{N_{k-1}}
u\|_{p_1}\|\triangle_{N_{k}} u\|_{p_1}^a . \label{inter8}
\end{align}
Inserting \eqref{inter8} into \eqref{inter-7}, we have
\begin{align} \|
u\|_{\dot{B}_{p, q}^{s}} & \lesssim  \sum\limits_{N_1\geq \ldots
\geq N_K}(N_1^{s }\ldots N_K^{s }\|\triangle_{N_1} u\|_p\ldots
\|\triangle_{N_K}
u\|_p)^{q^2} \nonumber\\
& \quad \times \Lambda(N_1,..., N_K) \| u\|_{\dot{B}_{p_1,
\infty}^{s_1}}^{q(1-q)\theta K} \| u\|_{\dot{B}_{p_0,
\infty}^{0}}^{(1-\theta)K q(1-q)}, \label{inter9}
\end{align}
where
\begin{align} \Lambda (N_1, \ldots N_K)& =
\Big(N_1^{-\frac{n}{p}+\frac{n}{p_1}-s_1+s}\ldots
N_{k-1}^{-\frac{n}{p}+\frac{n}{p_1}-s_1+s}N_k^{a(-\frac{n}{p}+\frac{n}{p_1}-s_1+s)}
\nonumber\\
& \quad \times
N_{k}^{(1-a)(-\frac{n}{p}+\frac{n}{p_0}+s)}N_{k+1}^{-\frac{n}{p}+\frac{n}{p_0}+s}
\ldots N_{K}^{-\frac{n}{p}+\frac{n}{p_0}+s} \Big )^{q(1-q)}.
\label{inter10}
\end{align}
By \eqref{inter9}, we have
\begin{align}
\| u\|_{\dot{B}_{p, q}^{s}} & \lesssim \sum\limits_{N_1\geq \ldots
\geq
N_K}\Lambda (N_1, \ldots N_K) \sum\limits_{i=1}^{K}(N_i^{s}\|\Delta_i u\|_p)^q \\
&\quad \times \| u\|_{\dot{B}_{p_1, \infty}^{s_1}}^{(1-q)\theta } \|
u\|_{\dot{B}_{p_0, \infty}^{0}}^{(1-\theta)(1-q)}. \label{inter11}
\end{align}
So, it suffices to prove
\begin{align}
\sum\limits_{N_1\geq \ldots \geq N_K}\Lambda (N_1, \ldots N_K)
\sum\limits_{i=1}^{K}(N_i^{s}\|\Delta_i u\|_p)^q \lesssim \|
u\|_{\dot{B}_{p, q}^{s}}^{q}. \label{inter12}
\end{align}
In fact, \eqref{inter10}--\eqref{inter12} imply the result. Finally,
we prove \eqref{inter12}. Applying the condition \eqref{inter-6},
we have
\begin{align}
 & \sum\limits_{N_1\geq \ldots \geq
N_K}\Lambda (N_1, \ldots N_K)(N_k^{s}\|\Delta_k u\|_p)^q \nonumber\\
&\lesssim \sum\limits_{N_{k-1}\geq N_k}\Big
(N_{k-1}^{(k-1)(s-s_1+\frac{n}{p_1}-\frac{n}{p})}N_k^{(K-k+1-a)(s
+\frac{n}{p_0}-\frac{n}{p})+a(s
-s_1+\frac{n}{p_1}-\frac{n}{p})}\Big)^{q(1-q)}
N^{sq}_k\|\Delta_k u\|_p^q \nonumber\\
&\lesssim \sum\limits_{N_{k-1}\geq
N_k}\left(\frac{N_{k-1}}{N_k}\right)^{(k-1)(s-s_1+\frac{n}{p_1}-\frac{n}{p})q(1-q)}
N^{sq}_k \|\Delta_k u\|_p^q \nonumber\\
&\lesssim  \| u\|_{\dot{B}_{p, q}^{s}}^{q}.
\end{align}

{\it Case} 2. We consider the case
\begin{align}
s_1-s +\frac{n}{p} -\frac{n}{p_1}<0, \quad  s+\frac{n}{p_0}-
\frac{n}{p}<0. \label{inter-17}
\end{align}
Substituting the summation $\sum_{N_1\ge...\ge N_K}$ by $\sum_{N_1\leq...\leq N_K}$
in \eqref{inter-7} and repeating the procedure as in Case 1, we can get the result, as desired.

Up to now, we have shown the results for the following two cases: (i) $s=(1-\theta)s_0 + \theta s_1$ and $p_0 =p_1$; (ii) $s <(1-\theta)s_0 + \theta s_1$ and $p \ge p_0 \vee p_1$.

\medskip

{\bf Step 2.} We consider the case $p< p_0\vee p_1$ and $s <(1-\theta)s_0 + \theta s_1$.   Due to $\theta \in (0,1)$ and $1/p \leq (1-\theta)/p_0 + \theta/p_1$, we see that $p_0\not= p_1$  and $p_0\wedge p_1 <p < p_0\vee p_1$.   Let $0< \varepsilon \ll 1$. In view of the result as in Step 1, we see that
\begin{align} \label{interp1}
\|f\|_{\dot B^s_{p,q}} \lesssim \|f\|^{1/2}_{\dot B^{s-\varepsilon}_{p,\infty}}\|f\|^{1/2}_{\dot B^{s+\varepsilon}_{p,\infty}}.
\end{align}
Since $s_0-n/p_0 \not= s_1-n/p_1$, we can assume that $s_0-n/p_0 < s_1-n/p_1$. It follows that $1/p-s/n \in (1/p_0-s_0/n, \; 1/p_1-s_1/n)$.
Hence, for sufficiently small $\varepsilon >0$,
$$
\frac{1}{p}- \frac{s\pm \varepsilon}{n} \in \left(\frac{1}{p_0}- \frac{s_0}{n}, \; \frac{1}{p_1}- \frac{s_1}{n} \right).
$$
It follows that there exist $\theta_\pm \in (0,1)$ satisfying
$$
\frac{1}{p}- \frac{s\pm \varepsilon}{n} = (1-\theta_\pm) \left(\frac{1}{p_0}- \frac{s_0}{n} \right) + \theta_\pm  \left(\frac{1}{p_1}- \frac{s_1}{n} \right).
$$
Due to $\lim_{\varepsilon \to 0} \theta_\pm = \theta$, we see that for sufficiently small $\varepsilon>0$,
$$
s\pm \varepsilon  \leq (1-\theta_\pm) s_0 +  \theta_\pm s_1.
$$
Therefore, by Theorem \ref{Besov-GN}, we have
\begin{align} \label{interp2}
& \|f\|_{\dot B^{s-\varepsilon}_{p,\infty}} \lesssim \|f\|^{1-\theta_-}_{\dot B^{s_0}_{p_0,\infty}} \|f\|^{\theta_-}_{\dot B^{s_1}_{p_1,\infty}},\\
& \|f\|_{\dot B^{s+\varepsilon}_{p,\infty}} \lesssim \|f\|^{1-\theta_+}_{\dot B^{s_0}_{p_0,\infty}} \|f\|^{\theta_+}_{\dot B^{s_1}_{p_1,\infty}}. \label{interp3}
\end{align}
We easily see that $\theta= (\theta_++\theta_-)/2$. Inserting \eqref{interp2} and \eqref{interp3} into \eqref{interp1}, we have the result, as desired.

(Necessity) First, we show the necessity for $s-n/p \neq s_0- n/p_0$. If not, then $s-n/p= s_0-
n/p_0=s_1 -n/p_1$. Let
\begin{align}
\hat{f} (\xi) = \sum^N_{j=100} 2^{(n/p-s)j} \varphi_j (\xi).
\end{align}
We see that $\|f\|_{B^s_{p,q}} \sim N^{1/q}$, $\|f\|_{B^{s}_{p,\infty}}
\sim 1$, which contradicts  \eqref{besovgn-3}.

\medskip

Next, we show the necessity of $p_0=p_1$ when $s=  (1-\theta) s_0 +  \theta  s_1. $  Assume for a contrary that $p_0 \neq p_1$.   By Theorem \ref{Besov-GN},  we have $1/q \leq (1-\theta)/\infty +  \theta /\infty=0$. This   contradicts the condition $q<\infty$. $\hfill \Box$

\medskip

 {\bf Proof of Theorem \ref{Trieb-GN}}  (Sufficiency) First, we consider the case $s< (1-\theta) s_0+
\theta s_1$. We can take sufficiently small $\varepsilon>0$ satisfying
$$
s\leq (1-\theta) s^*_0+ \theta s^*_1 , \quad s^*_0:=
s_0-\varepsilon, \; s^*_1:= s_1-\varepsilon.
$$
Since $\varepsilon \ll 1$, we can assume that
$$
\frac{1}{p^*_0}: =\frac{1}{p_0} -\frac{\varepsilon}{n}>0, \quad
\frac{1}{p^*_1}: =\frac{1}{p_1} -\frac{\varepsilon}{n}>0.
$$
Hence,
\begin{align}
\frac{n}{p}- s =
(1-\theta)\left(\frac{n}{p^*_0}- s^*_0 \right)
+\theta\left(\frac{n}{p^*_1}- s^*_1 \right), \label{}
\end{align}
which implies that
\begin{align}
\frac{1}{p}- \frac{1-\theta}{p^*_0}- \frac{\theta}{p^*_1}=
\frac{s}{n}- (1-\theta)\frac{s^*_0}{n}- \theta
\frac{s^*_1}{n}:=-\eta\leq 0.
\end{align}
Putting
\begin{align}
\frac{1}{p^*}= \frac{1}{p}+\eta, \quad s^*= s+n\eta,
\end{align}
we see that
\begin{align}
\frac{1}{p^*}= \frac{1-\theta}{p^*_0} +\frac{\theta}{p^*_1}, \quad
s^*= (1-\theta) s^*_0+ \theta s^*_1. \label{}
\end{align}
Using H\"older's inequality, in an analogous way as in Besov spaces, we have
$$
\|f\|_{\dot F^{s^*}_{p^*,q}} \lesssim \|f\|^{1-\theta}_{\dot
F^{s^*_0}_{p^*_0,q}} \|f\|^\theta_{\dot F^{s^*_1}_{p^*_1,q}}.
$$
Recalling the inclusions (see Triebel \cite{Tr})
$$
\dot F^{s_0}_{p_0,\infty}  \subset F^{s^*_0}_{p^*_0,q}, \quad \dot
F^{s_1}_{p_1,\infty}  \subset F^{s^*_1}_{p^*_1,q}
$$
we immediately get the conclusion.

Next, we consider the case $s= (1-\theta) s_0+
\theta s_1$ and $s_0 \neq s_1$. In this case we easily see that $1/p= (1-\theta)/p_0+ \theta/p_1$.  The result has been shown in \cite{Oru98} and \cite{BM01} and we omit the details of the proof.

(Necessity) It suffices to consider the necessity in the case  $s=(1-\theta)s_0+ \theta s_1$.  If not, then $s_0=s_1=s$.  Let $\rho_j$ be as in \eqref{conter}
and
\begin{align}
\hat{f} (\xi) = \sum^N_{j=100} 2^{-sj}\rho_j (\xi).
\end{align}
We easily see that
$$
\|f\|_{\dot F^s_{p,\infty}} =
\|\mathscr{F}^{-1} (\rho_0)\|_p \sim 1.
$$
But
$$
\|f\|_{\dot F^s_{p,q}} \sim N^{1/q},
$$
which contradicts the GN inequality.
$\hfill\Box$

\section{Corollaries of the GN inequalities}

In this section we give some corollaries of our main results. Noticing that $BMO =\dot F^0_{\infty,2} \subset \dot B^0_{\infty, \infty}$ and $\|\nabla^s u\|_{\dot B^0_{p,\infty}} \lesssim  \| \nabla^s u\|_{p}$,  we can deduce the following useful interpolation inequalities:
\begin{align}
& \|u\|_{L^{10} (\mathbb{R}^3)} \leq C\|u\|^{2/3}_{\dot
B^{-1/2}_{\infty, \infty}(\mathbb{R}^3) }\|u\|^{1/3}_{\dot
B^{1}_{10/3, 10/3}(\mathbb{R}^3)}, \label{Bouint}\\
& \|u\|_{L^4} \lesssim \|\nabla u\|^{1/2}_{L^2} \|u\|^{1/2}_{\dot
B^{-1}_{\infty, \infty}}, \label{MeyerRivint}\\
& \|\nabla u\|_{L^4} \lesssim \|\nabla^2 u\|^{1/2}_{L^2} \|u\|^{1/2}_{BMO}, \label{MeyerRivint2}\\
& \|u\|_{L^q} \lesssim \|\nabla u\|^{\theta}_{L^p} \|u\|^{1-\theta}_{\dot
B^{-\theta/(1-\theta)}_{\infty, \infty}}, \quad 1\le p<q<\infty, \theta=p/q. \label{Ledoux}\\
& \|\nabla^m u\|_{L^q} \lesssim \|\nabla^k u\|^{\theta}_{L^p} \|u\|^{1-\theta}_{BMO}, \quad 1\leq m<k, \ q=kp/m, \ \theta=m/k. \label{Strzelecki}
\end{align}
Following Bourgain \cite{Bourg1}, we can show \eqref{Bouint}, which is useful for the concentration phenomena for the solutions of nonlinear Schr\"odinger equations.  Meyer and Rivi\`{e}re \cite{MeRi} studied the partial regularity of solutions for the stationary Yang-Mills fields by using \eqref{MeyerRivint} and \eqref{MeyerRivint2}. \eqref{Ledoux} and \eqref{Strzelecki} are generalized versions of  \eqref{MeyerRivint} and \eqref{MeyerRivint2}, respectively (see Ledoux \cite{Led03}, Strzelecki \cite{Strz06}).   Machihara and Ozawa \cite{MaOz02} showed that
\begin{prop}
Let  $1\le   p_0 \vee p_1\leq p \leq  \infty$, $0<\theta<1$,
$  s_0,  s_1 \in \mathbb{R}$. Assume that
\begin{align}
\frac{n}{p}- s  & = (1-\theta)\left(\frac{n}{p_0}-
 s_0 \right)
+\theta\left(\frac{n}{p_1}- s_1 \right), \nonumber\\
& s_0 < \frac{n}{p_0} - \frac{n}{p }, \ \ s_1> \frac{n}{p_1} - \frac{n}{p }.
\end{align}
Then
\begin{align}
\| u\|_{\dot{B}_{p,1}^{0}}\lesssim  \| u\|_{\dot{B}_{p_0,
\infty}^{s_0}}^{1-\theta} \| u\|_{\dot{B}_{p_1, \infty}^{s_1}}^{
\theta}
\end{align}
\end{prop}

\medskip

\noindent Oru \cite{Oru98} obtained that (see also \cite{BM01})
\begin{prop}\label{Prop1.5.6.}
Let $0< p_0, p_1, p < \infty$, $0<r< \infty$,
$-\infty<s_0, s_1, s<\infty$, $0<\theta<1$ and
\begin{align}
\frac{1}{p}=\frac{1-\theta}{p_0}+\frac{\theta}{p_1}, \quad
s=(1-\theta) s_0+ \theta s_1, \ \ s_0\neq s_1.  \label{convex}
\end{align}
Then
\begin{align}
\|u\|_{\dot F^s_{p,r} (\mathbb{R}^n)} \leq C\|u\|^{1-\theta}_{\dot
F^{s_0}_{p_0, \infty}(\mathbb{R}^n) }\|u\|^{\theta}_{\dot
F^{s_1}_{p_1, \infty}(\mathbb{R}^n)}. \label{convex5}
\end{align}
\end{prop}

\medskip

\noindent The following interpolation inequality was shown in \cite{Wang5}.

\begin{prop}\label{Prop1.5.6.}
Let $0< p_0< p< \infty$, $0<r\leq \infty$,
$-\infty<s_1<s<s_0<\infty$, $0<\theta<1$ and
\begin{align}
\frac{1}{p}=\frac{\theta}{p_0}+\frac{1-\theta}{\infty}, \quad\quad
s=\theta s_0+ (1-\theta) s_1.  \label{convex1}
\end{align}
Then
\begin{align}
\|u\|_{\dot F^s_{p,r} (\mathbb{R}^n)} \leq C\|u\|^{1-\theta}_{\dot
B^{s_1}_{\infty, \infty}(\mathbb{R}^n) }\|u\|^{\theta}_{\dot
B^{s_0}_{p_0, p_0 }(\mathbb{R}^n)}. \label{convex2}
\end{align}
\end{prop}

\section{GN inequalities in nonhomogeneous spaces \label{nonhfs}}

We denote by $  H^s_p:= (I-\triangle)^{s/2} L^p$ the Bessel potential space, $  H^s=   H^s_2$.
Let $\psi$ be a smooth cut-off function supported in the ball $\{\xi:  |\xi| \leq 2\}$,  $\varphi= \psi(\cdot)- \psi(2\,\cdot)$.
We write $\psi_0: = \psi$ and
$ \psi_k (\xi)=\varphi (2^{-k} \xi), \  k\in \mathbb{N}.$
We see that
\begin{align}\label{eq1.5.2}
\sum^\infty_{k=0} \psi_k(\xi)=1,  \quad \xi\in \mathbb{R}^n.
\end{align}
We introduce the  dyadic decomposition operators
$\triangle_k=\mathscr{F}^{-1}\varphi_k \mathscr{F}, \quad  k\in \mathbb{Z}_+.$
Let
$-\infty<s<\infty, \;  1\leq p, q\leq \infty.$
The space $  B^s_{p,q}$ equipped with norm
\begin{align}\label{eq1.5.6}
 \|f\|_{
B^s_{p,q}}:=\bigg(\sum^\infty_{k=0} 2^{k s q}\|\triangle_k f
\|^q_p \bigg)^{1/q}
\end{align}
 is said to be a   Besov space. Let
\begin{align}\label{eq1.5.8}
-\infty <s<\infty, \quad 1\leq p< \infty, \quad 1\leq q\leq  \infty.
\end{align}
The space $  F^s_{p,q}$ equipped with norm
\begin{align}
 \label{eq1.5.10}
\|f\|_{F^s_{p,q}}:=\bigg\|\bigg(\sum^\infty_{k=0} 2^{k
s q} |\triangle_k f |^q\bigg)^{1/q} \bigg \|_p
\end{align}
is said to be a homogeneous Triebel-Lizorkin space. For Besov spaces and Triebel spaces, we have similar results as in Theorems \ref{Besov-GN}, \ref{Besov-GN1} and \ref{Trieb-GN}.  In this paper, we will use the following

\begin{thm} \label{nBesov-GN1}
Let $0<q <\infty$, $0< p, \; p_0, \; p_1\leq \infty$, $0<\theta<1$,
$-\infty<s, s_0,  s_1 <\infty$. Then the GN inequality of the following type
\begin{align}
\| u\|_{ {B}_{p, q}^{s}}\lesssim  \| u\|_{ {B}_{p_0,
\infty}^{s_0}}^{1-\theta} \| u\|_{ {B}_{p_1, \infty}^{s_1}}^{
\theta} \label{besovgn-3}
\end{align}
holds if
\begin{align}
\frac{n}{p}- s  & = (1-\theta)\left(\frac{n}{p_0}-
  s_0 \right)
+\theta\left(\frac{n}{p_1}- s_1 \right), \label{bespl0}\\
 s_0  -\frac{n}{p_0}& \neq s_1- \frac{n}{p_1},   \label{bespl}\\
 s & \leq  (1-\theta)s_0 + \theta s_1, \label{besovgn-4}\\
       p_0 & = p_1     \ \mbox{ if }   \  s= (1-\theta) s_0+ \theta s_1.     \label{besovgn-2022}
\end{align}
\end{thm}

\begin{prop} \label{Bessel-GN}
Let $1<p, p_0, p_1< \infty$, $s,\ s_1\in
 \mathbb{R}$,  $0  \leq \theta \leq 1$. Then the GN inequality of the following type
\begin{align}
\| u\|_{ {H}_{p}^{s}}\lesssim  \| u\|_{L^{p_0}}^{1-\theta} \|
u\|_{ {H}_{p_1}^{s_1}}^{\theta} \label{riesz-1}
\end{align}
holds if
\begin{align}
\frac{n}{p}- s = (1-\theta)\frac{n}{p_0}
+\theta\left(\frac{n}{p_1}- s_1 \right), \quad s \leq \theta
s_1. \label{riesz-2}
\end{align}
\end{prop}
The proofs of these results are the same as those in  homogeneous spaces and the details of the proofs are omitted.

\section{Concentration of solutions of NS equation}

The local well posedness in $L^n$ for the NS equation is well-known; cf. Kato \cite{Kato84}. Here we need the following result (see for instance \cite{KeKo09} in 3D and \cite{WHHG} in higher spatial dimensions).

\begin{thm}\label{Prop4.3.1}
Let $u_0 \in  L^n $ with ${\rm div} u_0 =0$. Then there exists a $T_m>0$ such that
the NS equation \eqref{eq4.0.1} has a unique solution $u$ satisfying
\begin{align}\label{eq4.3.13}
 u \in C ([0,T_m); \  L^n) \  \cap     L^{2+n}_{\rm loc}(0,T_m; \  L^{2+n}).
\end{align}
If $T_m<\infty$, then we have $\|u\|_{ L^{2+n}(0,T_m; \ L^{2+n})}=\infty$. Moreover, if $u_0\in L^2$, then
\begin{align}\label{eq4.3.14}
\frac{1}{2}\|u(t)\|^2_2 +   \int^t_0 \|\nabla u(s)\|^2_2
ds = \frac{1}{2} \|u_0\|^2_2, \  \ 0<t<T_m.
\end{align}
\end{thm}
We will sketch the proof of Theorem  \ref{Prop4.3.1} in the Appendix.
In the sequel,  we will write $\|u\|^2_2:= \sum^n_{i=1}\|u_i\|^2_2$, $\|\nabla u\|^2_2:= \sum^n_{i,j=1}\|\partial_{x_j} u_i\|^2_2$ for $u=(u_1,...,u_n)$.  We have the following
\begin{prop} \label{NS estimates}
Let $\sigma\ge 1$ and $u$ be the smooth solution of NS equation. Then we have
\begin{align}
\frac{1}{2+\sigma} \frac{d}{dt} \|u(t)\|^{2+\sigma}_{2+\sigma} & + \frac{1}{2} \int_{\mathbb{R}^n} (\nabla |u|^\sigma \cdot \nabla |u|^2)(x)dx
\nonumber \\
 & +  \int_{\mathbb{R}^n}   |u|^\sigma  |\nabla u|^2 (x )dx - \int_{\mathbb{R}^n}  (\nabla p \cdot |u|^\sigma u) (x)dx  =0.
\label{1}
\end{align}
\end{prop}

{\bf Proof.} The first equation in   \eqref{eq4.0.1} is multiplied by $|u|^\sigma u $, we have
\begin{align}\label{eq4.0.3}
  |u|^\sigma u \cdot \left(\partial_t u    - \Delta u  + \sum^n_{j=1}  u_j \partial_{x_j} u  +   \nabla p \right)=0.
\end{align}
We have
\begin{align}\label{eq4.0.4}
 \sum^n_{i=1} |u|^\sigma u_i  \partial_t u_i  = \frac{1}{2} |u|^\sigma  \partial_t |u|^2 =  \frac{1}{2+\sigma}  \partial_t |u|^{\sigma+2},
\end{align}
\begin{align}\label{eq4.0.5}
 |u|^\sigma u_i \Delta u_i   = \nabla (|u|^\sigma  u_i \nabla u_i) - \frac{1}{2} (\nabla |u|^\sigma \cdot \nabla u_i^2) - |u|^\sigma |\nabla u_i|^2.
   \end{align}
It follows that
\begin{align}\label{eq4.0.5}
\sum^n_{i=1} |u|^\sigma u_i \Delta u_i   = \sum^n_{i=1} \nabla (|u|^\sigma  u_i \nabla u_i) - \frac{1}{2} (\nabla |u|^\sigma \cdot \nabla |u|^2) - |u|^\sigma |\nabla u|^2.
\end{align}
Noticing that ${\rm div } u=0$, we have
\begin{align}\label{eq4.0.6}
\sum^n_{i=1} |u|^\sigma u_i \sum^n_{j=1} u_j \partial_j u_i   & = \frac{1}{2}\sum^n_{i,j=1}  |u|^\sigma  u_j  \partial_j u^2_i
   = \frac{1}{2}\sum^n_{j=1}  |u|^\sigma  u_j  \partial_j |u|^2  \nonumber \\
&   = \frac{1}{2+\sigma}\sum^n_{j=1} \partial_j (|u|^{\sigma+2}  u_j ) .
   \end{align}
We obtain that
\begin{align}\label{eq4.0.7}
 \frac{1}{2+\sigma}  \partial_t |u|^{\sigma+2} & -  \frac{1}{2} \nabla (|u|^\sigma    \nabla |u|^2) + \frac{1}{2} (\nabla |u|^\sigma \cdot \nabla |u|^2)  \nonumber\\
& + |u|^\sigma |\nabla u|^2  +   |u|^\sigma u   \nabla p  + \frac{1}{2+\sigma}  \nabla (|u|^{\sigma+2}  u ) =0.
   \end{align}
Integrating \eqref{eq4.0.7} over $\mathbb{R}^n$, we immediately obtain the result, as desired. $\hfill \Box$

Recall that by \eqref{eq4.0.1},
\begin{align}\label{eq4.3.15}
  -\Delta p= \sum^{n}_{i,j=1} \partial_{x_ix_j} (u_iu_j).
\end{align}
Let us denote
\begin{align}\label{eq4.3.16}
  E(u, v) = \sum^{n}_{i,j=1}   \mathscr{F}^{-1} |\xi|^{-2} \xi_i \xi_j \mathscr{F} (u_i v_j).
\end{align}
From the H\"ormander-Mikhlin multiplier theorem, we obtain that for any $p\in (1,\infty)$,
\begin{align}\label{eq4.3.17}
 \| E(u, v)\|_p  \lesssim  \sum^{n}_{i,j=1} \|  u_i v_j \|_p.
\end{align}

Putting $\sigma=n-2$  and integrating \eqref{1} over $[t_1,t_2]$, we have
\begin{align}
  \|u(t_2)\|^{n}_{n} &       +  2(n-2)  \int^{t_2}_{t_1} \int_{\mathbb{R}^n} \left|\nabla |u|^{n/2}\right |^2 dxdt
\nonumber \\
 &    +  n \int^{t_2}_{t_1} \int_{\mathbb{R}^n}   |u|^{n-2}  |\nabla u|^2  dx dt \le    \|u(t_1)\|^{n}_{n}  + n \int^{t_2}_{t_1} \int_{\mathbb{R}^n}  | p \nabla ( |u|^{n-2} u) | dx dt.
\label{1l}
\end{align}
Applying \eqref{eq4.3.14} and \eqref{eq4.3.17}, we obtain that
\begin{align}
 & \int^{t_2}_{t_1} \int_{\mathbb{R}^n}  \left| p \nabla ( |u|^{n-2}  u) \right | dx dt \nonumber\\
 & \lesssim   \int^{t_2}_{t_1}  \||u|^{(n-2)/2}  |\nabla u| \|_2 \||u|^{(n-2)/2} E(u,u)\|_2 dt \nonumber\\
 & \lesssim   \int^{t_2}_{t_1}  \||u|^{(n-2)/2}  |\nabla u| \|_2  \|u\|^{(n-2)/2}_{2+n} \| E(u,u)\|_{(n+2)/2} dt \nonumber\\
 & \lesssim   \frac{1}{100} \int^{t_2}_{t_1}  \||u|^{(n-2)/2}  |\nabla u| \|^2_2 dt + C_n \int^{t_2}_{t_1} \|u\|^{ 2+n}_{2+n}   dt.
\label{1ll}
\end{align}
Inserting the estimate as in \eqref{1ll} into \eqref{1l}, we have

\begin{lem}\label{lemns}
Let $u$ be the solution of
the NS equation \eqref{eq4.0.1} in $[0,T_m)$ and $t_1,t_2 \in [0,T_m)$. We have
\begin{align}
  \|u(t_2)\|^{n}_{n} &       +  2(n-2)  \int^{t_2}_{t_1} \int_{\mathbb{R}^n} \left|\nabla |u|^{n/2}\right |^2 dxdt
\nonumber \\
 &    + \frac{99n}{100}   \int^{t_2}_{t_1} \int_{\mathbb{R}^n}   |u|^{n-2}  |\nabla u|^2  dx dt \le    \|u(t_1)\|^{n}_{n}  + C_n \int^{t_2}_{t_1} \|u\|^{2+n}_{2+n}  dt.
\label{1l}
\end{align}
\end{lem}

\noindent{\bf Proof of Theorem \ref{NS concent}.}  By the local well posedness result  and Lemma \ref{lemns}, we see that if  $T_m<\infty$, then we have
\begin{align}
 \|u\|_{L^{2+n}_{x,t \in [0, T_m)}} =\infty.
\label{blowupsol}
\end{align}
In the following we give the details of the analysis to  $ \|u\|_{L^{2+n}_{x,t \in [0, T]}}$.  We have
\begin{align}
 \int^T_S \|u (t)\|^{2+n}_{2+n}dt  & =  \int^T_S \left\||u (t)|^{n/2} \right\|^{2(2+n)/n}_{2(2+n)/n} dt   \nonumber\\
  &  \le  \int^T_S \left\|P_{\le N} |u (t)|^{n/2} \right\|^{2(2+n)/n}_{2(2+n)/n} dt  +  \int^T_S \left\|P_{\ge N} |u (t)|^{n/2} \right\|^{2(2+n)/n}_{2(2+n)/n} dt.
\label{blowupsola}
\end{align}
For convenience, we write
$$
P_{\le N} f :=   \mathscr{F}^{-1} \psi(2^{-N} \xi) \mathscr{F},
$$
where $\psi$ is the smooth cut-off function supported in $\{\xi: |\xi| \le 2\}$ as before.
Using Bernstein's estimates and the $L^2$ bound of solutions, we see that
\begin{align}
 &   \int^T_S \left\|P_{\le N} |u (t)|^{n/2} \right\|^{2(2+n)/n}_{2(2+n)/n} dt \nonumber\\
&   \lesssim  (T-S) 2^{Nn^2/2}  \max_{t\in [S,T]} \left\|P_{\le N} |u (t)|^{n/2} \right\|_{4/n}^{2(n+2)/n}    \nonumber\\
&   \lesssim  (T-S) 2^{Nn^2/2}  \max_{t\in [S,T]} \left\|  u (t) \right\|^{(n+2) /2}_2   \lesssim  (T-S) 2^{Nn^2/2}.
\label{blowupanal1}
\end{align}
Let $T_k \nearrow T_m$, we  see that  $\|u\|_{L^{2+n}_{x,t \in [0, T_k]}}  \nearrow \infty$. We can assume, by passing to a subsequence of $\{T_k\}$ that  \begin{align}
\|u\|_{L^{2+n}_{x,t \in [T_{k-1}, T_k]}}  \ge \|u\|_{L^{2+n}_{x,t \in [0, T_{k-1}]}}. \label{mon}
\end{align}
Let $N_k \nearrow \infty$ satisfy
\begin{align}
c \|u\|^{2+n}_{L^{2+n}_{x,t \in [T_{k-1}, T_k]}} \le  C (T_k-T_{k-1}) 2^{N_k n^2/2}  \le  \frac{1}{2} \|u\|^{2+n}_{L^{2+n}_{x,t \in [T_{k-1}, T_k]}}.
\label{blowupanal2}
\end{align}
We have
\begin{align}
    \int^{T_k}_{T_{k-1}} \left\|P_{\le N_k} |u (t)|^{n/2} \right\|^{2(2+n)/n}_{2(2+n)/n} dt
 \le  \frac{1}{2} \|u\|^{2+n}_{L^{2+n}_{x,t \in [T_{k-1}, T_k]}}.
\label{blowupanal3}
\end{align}
It follows from \eqref{blowupsola}, \eqref{mon} and \eqref{blowupanal3} that
\begin{align}
 \frac{1}{4} \|u\|^{2+n}_{L^{2+n}_{x,t \in [0, T_k]}} \le    \int^{T_k}_{T_{k-1}} \left\|P_{\ge N_k} |u (t)|^{n/2} \right\|^{2(2+n)/n}_{2(2+n)/n} dt.
\label{blowupanal4}
\end{align}
In view of the fractional GN inequality, we have
\begin{align}
  \left\|v  \right\|_{L^{2(2+n)/n}_{x,t\in [{T_{k-1}} ,T_k]}} \lesssim \|v\|^{2/(n+2)}_{L^\infty (T_{k-1}, T_k; \dot B^{-n/2}_{\infty, \infty})} \|\nabla v\|^{n/(n+2)}_{L^{2}_{x,t\in [T_{k-1}, T_k]}}.
\label{blowupanal5}
\end{align}
Taking $v= P_{\ge N_k} |u|^{n/2}$, by \eqref{blowupanal4} and \eqref{blowupanal5} we have
\begin{align}
 \frac{1}{4} \|u\|^{2+n}_{L^{2+n}_{x,t \in [0, T_k]}}   \lesssim \| P_{\ge N_k} |u|^{n/2} \|^{4/n}_{L^\infty (T_{k-1}, T_k; \dot B^{-n/2}_{\infty, \infty})} \|\nabla P_{\ge N_k} |u|^{n/2} \|^{2}_{L^{2}_{x,t\in [T_{k-1},T_k]}}.
\label{blowupanal6}
\end{align}
By Lemma \ref{lemns}, we see that
\begin{align}
 \|\nabla P_{\ge N_k} |u|^{n/2} \|^{2}_{L^{2}_{x,t\in [0,T_k]}} \lesssim  \|u_0\|^n_n + \|u\|^{2+n}_{L^{2+n}_{x,t \in [0, T_k]}}.
\label{blowupanal7}
\end{align}
Hence, it follows from  \eqref{blowupanal6} and \eqref{blowupanal7} that
 \begin{align}
 \| P_{\ge N_k} |u|^{n/2} \|_{L^\infty (T_{k-1}, T_k; \dot B^{-n/2}_{\infty, \infty})} \gtrsim 1.
\label{blowupanal8}
\end{align}
We remark that the constant in the right hand side of \eqref{blowupanal8}  only depends on $n$ and $\|u_0\|_n$. So, there exist $x_k \in \mathbb{R}^n$, $t_k \in [T_{k-1}, T_k]$ and $j_k \ge N_k$ such that
 \begin{align}
  2^{-nj_k/2} | (\Delta_{j_k} |u|^{n/2})(x_k, t_k)  | \gtrsim 1.
\label{blowupanal9}
\end{align}
Let $\psi$ be as in \eqref{eq1.5.2}, $0<\varepsilon \ll 1$. It follows that
\begin{align}
1 & \lesssim   2^{n j_k/2}  \left| \int (\mathscr{F}^{-1} \psi)(2^{j_k}(x_k-y)) |u(t_k, y) |^{n/2} dy \right| \nonumber\\
& \lesssim   2^{n j_k/2}  \left| \left(\int_{|y-x_k|\le 2^{(\varepsilon-1)j_k}} + \int_{|y-x_k| > 2^{(\varepsilon-1)j_k}} \right) (\mathscr{F}^{-1} \psi)(2^{j_k}(x_k-y)) |u(t_k, y) |^{n/2} dy \right| \nonumber\\
&:= I+II. \label{blowupanal10}
\end{align}
By H\"older's inequality, we have
\begin{align}
 II   & \lesssim   2^{n j_k/2}  \|(\mathscr{F}^{-1} \psi)(2^{j_k}\, \cdot) \|_{L^{4/(4-n)}(|\cdot|> 2^{(\varepsilon-1)j_k}) }\|u(t_k, \cdot)\|^{n/2}_{2} \nonumber\\
& \lesssim   2^{n j_k/2}  \|(\mathscr{F}^{-1} \psi)(2^{j_k}\, \cdot) \|_{L^{4/(4-n)}(|\cdot|> 2^{(\varepsilon-1)j_k})} \nonumber\\
&  \lesssim   2^{n j_k(n/4- 1/2) }   \|(\mathscr{F}^{-1} \psi)  \|_{L^{4/(4-n)}(|\cdot|> 2^{ \varepsilon j_k})}. \label{blowupanal11}
\end{align}
Since $\psi$ is a Schwartz function, for fixed $\varepsilon>0$, we have
$$
II \ll 1/2, \quad if \quad k \gg 1.
$$
Hence,
$$
1/2 \lesssim I \lesssim \|u\|^{n/2}_{L^n(|\cdot-x_k| \le  2^{ (\varepsilon-1) j_k})} \|\mathscr{F}^{-1} \psi\|_2  \lesssim  \|u\|^{n/2}_{L^n(|\cdot-x_k|\le 2^{ (\varepsilon-1) j_k})}.
$$
By \eqref{blowupanal2}, we see that $2^{(\varepsilon-1)j_k} \lesssim (T_m -T_{k-1})^{\delta}$ for any $\delta< 2/n^2$.  $\hfill\Box$

\section{Proof of Theorem \ref{exist-minimizer}}

Let $q=n/(n-\beta)$. First, we consider the case $m^2=0$.  We divide the proof into the following five steps.

{\bf Step 1.} We show that $M_c>-\infty$. Applying  Hardy-Littlewood-Sobolev's inequality, we have
\begin{align}
 & \int \int G( u_1(x),...,u_L(x) ) V(|x-y|) G( u_1(y) ,..., u_L(y)  ) dxdy  \nonumber\\
 & \lesssim   (\|u\|^2_{2(2q)'}+ \|u\|^\mu_{(2q)'\mu})^2,
\end{align}
where $(2q)'$ is the dual exponent to $2q$.  In view of the fractional Gagliardo-Nirenberg inequality, we have
\begin{align}
  \|u\|_{2(2q)'} & \lesssim \|u\|^{1-\theta_2}_{2} \|u\|^{\theta_2}_{\dot H^s}, \\
   \|u\|_{\mu (2q)'} & \lesssim \|u\|^{1-\theta_\mu}_{2} \|u\|^{\theta_\mu}_{\dot H^s},
\end{align}
where
$$
\frac{s \theta_\lambda}{n}=  \frac{1}{2} -\frac{1}{\lambda (2q)'}.
$$
We consider the following two cases. First, if $\mu < 2+2s/n -1/q $,  we easily see  that $2\theta_2, \mu \theta_\mu <1$. It follows from $u\in S_c$ that
\begin{align} \label{nlc}
 & \int \int G( u_1(x) ,..., u_L(x) ) V(|x-y|) G( u_1(y),..., u_L(y)  ) dxdy  \nonumber\\
 & \lesssim \|u\|^{2\theta_\mu}_{\dot H^s}   + \|u\|^{2\theta_2}_{\dot H^s} \lesssim 1+ \varepsilon \|u\|^{2}_{\dot H^s}
\end{align}
for some sufficiently small $\varepsilon>0$. Next, if  $\mu = 2+2s/n -1/q $, applying the condition $u\in S_c$ and $c_1,..., c_L$ are sufficiently small, we see that \eqref{nlc} also holds. So, we have shown that
\begin{align} \label{nlceu}
E(u) \ge  \left(\frac{1}{2}     - C \varepsilon \right) \|u\|^{2}_{\dot H^s} -C.
\end{align}
Therefore, we have $M_c> -\infty$ and all of the minimizing sequence of \eqref{ques-minizer} are bounded in $(H^s)^L $.

{\bf Step 2.} We show the existence of the Schwarz symmetric (=radial and  radially decreasing) sequence. Let $u^*$ be the monotone rearrangement of $u$. By the super-modularity of $G$ (see Proposition 3.13 of \cite{HaKr}) and Theorem 1.2 in \cite{BuHa},
\begin{align} \label{nlcsuper}
 & \int \int G( u_1(x) ,..., u_L(x)  ) V(|x-y|) G( u_1(y) ,..., u_L(y)) dxdy  \nonumber\\
 & \leq \int \int G( u^*_1(x) ,..., u^*_L(x)) V(|x-y|) G( u^*_1(y),..., u^*_L(y)) dxdy.
\end{align}
On the other hand, we know thanks to  (cf. Appendix of \cite{AlLi})
\begin{align} \label{normcom}
\|u^*\|_{\dot H^s}  \leq \|u \|_{\dot H^s}.
\end{align}
It follows that $E(u^*) \leq E(u)$. Hence, we obtain the existence of the Schwarz minimizing sequence. So, it suffices to consider the Schwarz minimizing sequence below.

{\bf Step 3.} We show the lower semi-continuity of $E(\cdot)$ under the Schwarz minimizing sequence. Let $u_k=(u_{k,1},..., u_{k,L})$ be a Schwarz symmetric minimizing sequence. We show that if $u_k$ weakly converges to $u$ in $(H^s)^L $, then
\begin{align} \label{lowercon}
E(u) \leq \lim\inf E(u_k).
\end{align}
Since the minimizing sequence in $(H^s)^L $ is bounded, we see that there exists a subsequence, which is still written by $u_k$ such that $u_k$ weakly converges to $u=(u_1,..., u_L) $ in $(\dot H^s)^L$. It follows that
\begin{align} \label{lowercont2}
\|u \|^2_{\dot H^s} \leq \lim \inf_{k\to \infty} \|u_k\|^2_{\dot H^s}.
\end{align}
In the following we show that
\begin{align} \label{limG}
   \lim_{k\to \infty}  \int \int G(u_k (x) ) V(|x-y|) G(u_k(y) ) dxdy
 =  \int \int G( u (x)) V(|x-y|) G( u (y)) dxdy.
\end{align}
The sequence $u_k$ is bounded in $(H^s)^L $, so is in $L^{(2q)'\mu} \cap L^{2(2q)' }$. Since $u_k$ is a symmetric sequence, we can certainly find a subsequence of $u_k$ still written by $u_k$ such that $u_{k} \to u $ and $|u_{k,j}| \leq a_j$ for some $a_j \in L^{(2q)'\mu} \cap L^{2(2q)' }$. By the continuity of $G$ we have
$$
G(u_k(x))V(|x-y|) G(u_k(y)) \to  G(u (x))V(|x-y|) G(u (y)), \ \ k \to \infty
$$
for all $x,y \in \mathbb{R}^n$. On the other hand, since $G$ is non-decreasing with respect to all variables, we have from condition (G1) that
\begin{align} \label{Gineq1}
 G(u_k(x))V(|x-y|) G(u_k(y))
 & \leq G(a (x))V(|x-y|) G(a (y)) \nonumber\\
 & \lesssim    (|a (x)|^2+ |a(x)|^\mu) V(|x-y|) (|a (y)|^2+ |a(y)|^\mu).
 \end{align}
It follows that
\begin{align} \label{Gineq2}
& \int\int G(u_k(x))V(|x-y|) G(u_k(y)) dxdy \nonumber\\
  & \lesssim  \int\int  (|a (x)|^2+ |a(x)|^\mu) V(|x-y|) (|a (y)|^2+ |a(y)|^\mu) \nonumber\\
  & \lesssim (\|a\|^2_{2(2q)'} + \|a\|^\mu_{\mu(2q)'}) < \infty.
 \end{align}
In view of the dominated convergence theorem, we immediately have \eqref{limG}.

{\bf Step 4.} We show the strict negativity of $M_c$. Let $\varphi: \mathbb{R}^n \to (0,1)$  be a Schwarz radial function satisfying $\|\varphi \|_2= 1$. Taking $\varphi_i=c_i \varphi$, $i=1,...,L $ and  $\Phi_\lambda = \lambda^{n/2} \Phi (\lambda \cdot) = \lambda^{n/2} ( \varphi_1 (\lambda \cdot),... , \varphi_L (\lambda \cdot) )$.   Clearly, we have $ \Phi_{\lambda} \in S_c$. For convenience, we write $\alpha=\alpha_1+...+\alpha_L.$ We have from the second growth condition in (G1) that for $0<\lambda \ll 1$,
\begin{align} \label{Gineq2}
E(\Phi_{\lambda}) &  = \frac{1}{2} \|\Phi_{\lambda}\|^2_{\dot H^s}  - \int\int G(\Phi_{\lambda} (x)) V(|x-y|)  G(\Phi_{\lambda}(y)) dxdy \nonumber\\
  & =  \frac{1}{2} \lambda^{2s} \|\Phi_1 \|^2_{\dot H^s}  - \lambda^{-2n} \int\int G(\lambda^{n/2} \Phi  (x)) V(|x-y|/\lambda)  G(\lambda^{n/2} \Phi (y)) dxdy \nonumber\\
  & \leq   \frac{1}{2} \lambda^{2s} \|\Phi_1 \|^2_{\dot H^s}  - C \lambda^{- n + \alpha n - \beta}   \int\int
  \varphi   (x)^{\alpha} V(|x-y|)  \varphi  (y)^{\alpha}  dxdy \nonumber\\
  & \leq     \lambda^{2s} (C_1  - C_2 \lambda^{- n + \alpha n - \beta -2s})
 \end{align}
for some $C_1, C_2>0$. Noticing that $ n+\beta-n\alpha+2s >0$ and taking $0<\lambda \ll 1$, we immediately have $E(\Phi_{\lambda})<0$. It follows that $M_c<0$.

{\bf Step 5.} We show that $M_c$ is achieved. Notice that $M_c = E(u)$. It suffices to show that $\|u_i\|^2_2 =c_i$.   The strict negativity of $M_c$ and condition (G2) imply that $u_i\neq 0$ for all $i=1,..,L$. Let $t_i= c_i/\|u_i\|^2_2 $, $i=1,...,L$. We have $t_i\ge 1$ and $(t_1u_1,..., t_Lu_L) \in S_c$. Therefore,
$$
M_c \leq E(t_1u_1,..., t_L u_L) \leq t^2_{\max} E(u) = t^2_{\max}  M_c.
$$
Since $M_c<0$, we immediately have $t_{\max}=1$ and so, $t_1=...=t_L =1$. It follows that $u$ is a minimizer.

\medskip

 Next, we consider the case $m^2>0$. Since $\|u\|_{ \dot H^s} \le \|u\|_{H^s}$, $s>0$,  we see that the proof in Steps 1--3 and 5 holds if we substitute $\dot H^s$ by $H^s$.  Moreover,   noticing that $\|u_\lambda\|^2_{H^s} \lesssim \|u\|^2_2 + \lambda^{2s}\|u\|^2_{\dot H^s}$, we see that the result in Step 4 is also true.  $\hfill\Box$

 \medskip

 In the proof of Step 4, we easily see that for the single power case $G(|u|^2)= u^\alpha$ with $2s+n+\beta-\alpha <0$, $E(\Phi_\lambda)\to -\infty$ as $\lambda \to \infty$. Moreover, taking $\alpha=2$, we see that the condition $s\ge (n-\beta)/2$ is also necessary.

\section{Proof of Theorem \ref{exist-minimizer2}}

(Necessity) Put $u_\lambda = \lambda^{n/2} u(\lambda \, \cdot)$, $s=(n-\beta)/2$. For any $\phi \in (H^{s})^L$, we write
\begin{align}
I^{(n)}_{c, \beta}  (\phi  )  =  \frac{1}{2}\|\phi\|^2_{\dot H^s}-  \Upsilon_\beta (\phi). \label{I}
\end{align}
we have
\begin{align}
I^{(n)}_{c, \beta}  (\phi_\lambda )  = \lambda^{n-\beta}  \left( \frac{1}{2}\|\phi\|^2_{\dot H^s}-  \Upsilon_\beta (\phi)  \right). \label{ChPe2lambda}
\end{align}
By \eqref{C*}
$$
\Upsilon_\beta (u)= \int_{\mathbb{R}^{2n}}  \frac{|u  (x)|^2  |u (y)|^2 }{|x-y|^{n-\beta} }dxdy \le  C^* c \|u\|^2_{\dot H^s}.
$$
If $  C^*c < 1/2$, then
\begin{align}
   \left(\frac{1}{2} - C^*c \right) \|u\|^2_{\dot H^s}  \le I^{(n)}_{c, \beta}  (u) \leq \frac{1}{2}    \|u\|^2_{\dot H^s}. \label{ChPe2ineq}
\end{align}
It follows that $M^{(n)}_{c, \beta} \ge 0.$ On the other hand, noticing that $\|\phi_\lambda\|_2 = \|\phi\|_2$, we see that
$$
M^{(n)}_{c, \beta}  \le \inf \{ I^{(n)}_{c, \beta} (\phi_\lambda ): \|\phi\|^2_2=c\} \leq \frac{ \lambda^{n-\beta}}{2}    \|\phi\|^2_{\dot H^s}
$$
holds for all $\phi \in H^s $  with $\|\phi\|^2_2=c$. Hence, $M^{(n)}_{c, \beta}=0$. For any minimizing sequence $u_k$, we have $I^{(n)}_{c, \beta}  (u_k) \sim   \|u_k\|^2_{\dot H^s} \to 0$. It follows that $u_k \to 0 $ in $(\dot H^s)^L$. But this contradicts the fact $\|u_k\|^2_2 =c$.

If $C^* c >1/2$, we have $(C^*-\varepsilon) c >1/2$ for sufficiently small $\varepsilon>0$. By the definition of $C^*$ we can choose some $\phi \in (H^s)^L $ such that
$$
\Upsilon_\beta (\phi) \ge  (C^*-\varepsilon) \|\phi\|^2_2 \|\phi\|^2_{\dot H^s}.
$$
However,
\begin{align}
I^{(n)}_{c, \beta}  (\phi_\lambda ) \le  \lambda^{n-\beta}  \left( \frac{1}{2}-   (C^*-\varepsilon ) c \right) \|\phi\|^2_{\dot H^s}. \label{ChPe2i}
\end{align}
Taking $\lambda \to \infty$, we immediately have $M^{(n)}_{c, \beta} =-\infty$. \\

\noindent (Sufficiency)  First, we show that $M^{(n)}_{c, \beta} =0$. Since $C^* c=1/2$, we have
$$
\Upsilon_\beta (u ) \leq \frac{1}{2} \|u\|^2_{\dot H^s}.
$$
It follows that $M^{(n)}_{c, \beta} \geq 0$. On the other hand, for any $\varepsilon >0$, we find some $\phi \in (\dot H^s)^L$ satisfying
$$
\Upsilon_\beta (\phi) \geq \frac{1-\varepsilon}{2} \|\phi\|^2_{\dot H^s}.
$$
For $s=(n-2)/2$, the above inequality is invariant under the scaling $\phi \mapsto \lambda^{n/2} \phi (\lambda \, \cdot)$, which implies that we can assume that $\|\phi\|_{\dot H^s} =1$. It follows that $I^{(n)}_{c, \beta}  (\phi) \leq \varepsilon$.  Hence $M^{(n)}_{c, \beta} =0$.

Now, let $u_k$ be a sequence verifying
\begin{align}  \label{miniseq}
\frac{\Upsilon_\beta (u_k) }{\|u_k\|^2_2 \|u_k\|^2_{\dot H^s}} \ge  C^* \left(1- \frac{1}{k} \right) .
\end{align}
Let $u^*_k$ be the rearrangement of $u_k$. Using the fact that
$$
\Upsilon_\beta (u_k) \le \Upsilon_\beta (u^*_k), \   \|u^*\|_{\dot H^s}  \leq \|u \|_{\dot H^s}, \  \|u^*\|_{2}  = \|u \|_{2},
$$
we see that \eqref{miniseq} also holds if    $u_k$ is replaced by $u^*_k$, i.e.,
\begin{align}  \label{miniseqsym}
\frac{\Upsilon_\beta (u^*_k) }{\|u^*_k\|^2_2 \|u^*_k\|^2_{\dot H^s}} \ge  C^* \left(1- \frac{1}{k} \right).
\end{align}
One can find $\lambda_k >0$ such that
$\|\lambda^{n/2}_k u^*_k (\lambda_k \, \cdot)\|_{\dot H^s} =1$. Since \eqref{miniseqsym} is invariant under the scaling $u^*_k \mapsto \lambda^{n/2} u^*_k (\lambda \, \cdot)$, we see that for $v_k =\lambda^{n/2}_k u^*_k (\lambda_k \, \cdot) $,
\begin{align}  \label{miniseqsymdec}
\frac{\Upsilon_\beta (v_k) }{\|v_k\|^2_2 \|v_k\|^2_{\dot H^s}} \ge  C^* \left(1- \frac{1}{k} \right)
\end{align}
and $\|v_k\|^2_2 =c$, $\|v_k\|_{\dot H^s}=1$. The inequality \eqref{miniseqsymdec} also implies that $I^{(n)}_{c, \beta}  (v_k) \leq 1/2k \to 0$.  It follows that $v_k$ is a radial and radially decreasing minimizing sequence. In view of $\|v_k\|^2_{H^s} \leq 1+c$ we see that $v_k$ has a subsequence which is still written by $v_k$ such that $v_k$  converges to $v$ with respect to the weak topology in $(H^s)^L $. On the other hand, the embedding $H^s \subset L^q$ with $s=(n-\beta)/2$, $2<q<2n/\beta$  is compact for the class of radial functions, we see that $v_k$ strongly converges to $v$ (up to a subsequence) in $(L^q)^L$ for all $2<q<2n/\beta$.  By \eqref{miniseqsym} and Theorem \ref{Besov-GN1}, we have for $k\geq 2$,
\begin{align}  \label{miniseqsyma}
1/4 \leq  \Upsilon_\beta (v_k)  \leq C \||v_k|^2\|^2_{2n/(n+\beta)} \le C \|v_k\|^2_{B^{s}_{2, \infty}} \|v_k\|^2_{B^{-n/2}_{\infty, \infty}}.
\end{align}
It follows that $ \|v_k\|_{B^{-n/2}_{\infty, \infty}} \geq c_0$, where $c_0:= 1/2\sqrt{C} $ is independent of $k$. Let $v_k = (v_k^1,...,v_k^L)$. It is easy to see that there exist $i\in \{1,2,...,L\}$ and a subsequence of $v^i_k$ which is still written by $v^i_k$ verifying $\|v^i_k\|_{B^{-n/2}_{\infty, \infty}} \geq c_0/L$.  From the definition of $B^a_{\infty, \infty}$ we can choose $j_k \in \mathbb{Z}_+$ and $x_k \in \mathbb{R}^n$,
\begin{align}  \label{c0}
c_0 /2L \leq   2^{-n j_k/2}  |(\triangle_{j_k} v^i_k) (x_k)|.
\end{align}
Denoting
$$
\mathbb{A}(j_k):= \{x: \ |x_k -x| \le A 2^{-j_k}\}, \ \
$$
\begin{align}  \label{c0est}
c_0 /2m & \leq   2^{-n j_k/2}  |(\mathscr{F}^{-1}\varphi_{j_k}) * v_k^i  (x_k)| \nonumber\\
&  =  2^{ n j_k/2}  \int_{\mathbb{R}^n}  (\mathscr{F}^{-1}\varphi_{j_k}) (2^{j_k}(x_k -z)) v_k^i(z) dz \nonumber\\
&  =  2^{ n j_k/2}  \left(\int_{\mathbb{A}(j_k)} + \int_{\mathbb{R}^n \setminus \mathbb{A}(j_k)}\right) (\mathscr{F}^{-1}\varphi_{j_k}) (2^{j_k}(x_k -z)) v_k^i(z) dz \nonumber\\
& :=I+II.
\end{align}
Taking $A:= A(\varphi,c) \gg 1$, we see that
$$
II \leq \|v_k^i\|_2 \|\mathscr{F}^{-1}\varphi\|_{L^2(|\cdot-x_k| \ge A )} \le c_0/4.
$$
By H\"older's inequality, we have
$$
I  \leq C\|v_k^i\|_{L^2(|\cdot-x_k| \le  A 2^{-j_k} )} \leq C\|v_k^i\|_{L^2(|\cdot-x_k| \le  A)}.
$$
We have
$$
 \|v_k^i\|_{L^2(|\cdot-x_k| \le  A)} \ge c_0/4C.
$$
Since $v_k^i$ is radial, we have $|x_k| \le X_0:= X_0(c_0, C, A)$. Indeed, in the opposite case we will have $\|v_k^i\|^2_2 > c$  if $|x_k|\gg 1$. So, we further have
$$
 \|v_k^i\|_{L^2(|\cdot| \le X_0+ A)} \ge c_0/4C.
$$
By H\"older's inequality,
$$
 \|v_k^i\|_{L^q(|\cdot| \le X_0+ A)} \ge \tilde{c}_0, \ \ \tilde {c}_0:= \tilde{c}_0(A, X_0, c_0).
$$
Since $v_k \to v$ in $(L^q)^L$, $2<q<2n/\beta$, we immediately have $v\neq 0$. Using the same way as in the proof of Theorem \ref{exist-minimizer}, we can get that
$$
0 \le I^{(n)}_{c, \beta}  (v) \le I^{(n)}_{c, \beta}  (v_k) \to 0.
$$
It follows that $I^{(n)}_{c, \beta}  (v)=0$. To finish the proof, it suffices to show that $\|v\|^2_2 =c$. If not, then we have $\|v\|^2_2 <c$.  Putting $\tilde{v} = \sqrt{c} v/\|v\|_2$, we have
\begin{align}
I^{(n)}_{c, \beta}  (\tilde{v} ) & =  \frac{c}{\|v\|^2_2}  \left( \frac{1}{2}\|v\|^2_{\dot H^s}-  \frac{c}{\|v\|^2_2}  \Upsilon (v)  \right) \nonumber\\
& = \frac{c}{\|v\|^2_2} I^{(n)}_{c, \beta}  (v) - \left(   \frac{c}{\|v\|^2_2} -1 \right) \Upsilon_\beta (v) <0,
\label{tildv}
\end{align}
which contradicts the fact that $I^{(n)}_{c, \beta}  (u) \ge 0$ for all $u\in (H^s)^L $.

\section{Proof of Theorem \ref{exist-minimizer4}}

We consider the variational problem
\begin{align}
M^{(n)}_{c, \beta,m}    & = \inf\{I^{(n)}_{c, \beta,m}  (u): \, u\in (H^s)^L , \ \|u\|^2_2=c>0 \},\\
I^{(n)}_{c, \beta,m}  (u)  & =    \frac{1}{2} \int (m^2+ |\xi|^2)^s |\widehat{u}(\xi)|^2 d\xi-  \Upsilon_\beta (u).
 \end{align}

\begin{lem} \label{Mcinfty}
Let $s=(n-\beta)/2$. If $C^*c >1/2$, then $M^{(n)}_{c, \beta,m} = - \infty.$
\end{lem}
{\bf Proof.} By Theorem \ref{exist-minimizer2}, there exists $\phi \in (H^s)^L $ with $\|\phi\|^2_2 =c$ satisfying
$$
\Upsilon_\beta (\phi)= C^* c \|\phi\|^2_{\dot H^s}.
$$
It follows that
\begin{align}
I^{(n)}_{c, \beta,m}  (\phi_\lambda)  & =    \frac{1}{2} \int (m^2+ |\lambda \xi|^2)^s |\widehat{\phi}(\xi)|^2 d\xi-  \lambda^{n-\beta}\Upsilon_\beta (\phi) \nonumber\\
& =    \frac{1}{2} \int (m^2+ |\lambda \xi|^2)^s |\widehat{\phi}(\xi)|^2 d\xi-  \lambda^{n-\beta} C^* c \|\phi\|^2_{\dot H^s}.   \label{I1}
 \end{align}
If $s\leq 1$, then
\begin{align}
I^{(n)}_{c, \beta,m} (\phi_\lambda)  \le   \frac{1}{2} m^{2s}  +      \lambda^{n-\beta}  \left( \frac{1}{2}- C^* c\right)  \|\phi\|^2_{\dot H^s}.
 \end{align}
Taking $\lambda\to \infty$, we immediately have $M^{(n)}_{c, \beta,m} =-\infty$.

Next, we consider the case $s>1$. Denote
$$
\mathbb{A} = \{\xi: \ \lambda |\xi| > m / \varepsilon\}.
$$
We have
\begin{align}
    \frac{1}{2} \int_{\mathbb{A}} (m^2+ |\lambda \xi|^2)^s |\widehat{\phi}(\xi)|^2 d\xi &  \leq  \frac{1}{2} \lambda^{2s} (1+ \varepsilon^{2})^s \int_{\mathbb{A}} |\xi|^{2s} |\widehat{\phi}(\xi)|^2 d\xi \nonumber\\
    &  \leq  \frac{1}{2} \lambda^{2s} (1+ \varepsilon^{2})^s  \|\phi\|^2_{\dot H^s}. \label{I2}
    \end{align}
On the other hand,
\begin{align}
    \frac{1}{2} \int_{\mathbb{R}^n \setminus \mathbb{A}} (m^2+ |\lambda \xi|^2)^s |\widehat{\phi}(\xi)|^2 d\xi &  \leq  \frac{1}{2} m^{2s} (1+  1/\varepsilon^2)^s \int_{\mathbb{R}^n \setminus \mathbb{A}}   |\widehat{\phi}(\xi)|^2 d\xi \nonumber\\
    &  \leq  \frac{1}{2} m^{2s} (1+ 1/\varepsilon^2 )^s  \|\phi\|^2_{2}. \label{I3}
    \end{align}
Collecting the estimates as in \eqref{I1}, \eqref{I2} and \eqref{I3}, we have
\begin{align}
I^{(n)}_{c, \beta,m}  (\phi_\lambda)  &  \le C_\varepsilon -    \lambda^{n-\beta}\left( C^* c - \frac{1}{2} (1+  \varepsilon^2)^s \right) \|\phi\|^2_{\dot H^s}.   \label{I4}
 \end{align}
By taking $\varepsilon>0$ small enough and $\lambda \to \infty$, we immediately have $M^{(n)}_{c, \beta,m} =-\infty$. $\hfill\Box$

\begin{lem}
Let $s=(n-\beta)/2 \ge 1$. If $C^*c \le 1/2$, then $M^{(n)}_{c, \beta,m} = cm^{2s}/2.$
\end{lem}
{\bf Proof.}   If $s\ge 1$, then we have
$$
(m^2+ |\xi|^2)^s \geq m^{2s} + |\xi|^{2s}.
$$
It follows that
$$
I^{(n)}_{c, \beta,m}  (\phi)     \geq \frac{1}{2} m^{2s}  \|\phi\|^2_{2} +   \frac{1}{2}   \|\phi\|^2_{\dot H^s} -\Upsilon_\beta (\phi).
$$
If $C^* c \le 1/2$ and $\|\phi\|^2_2 = c$, then we have
$$
\Upsilon_\beta (\phi) \le C^* \|\phi\|^2_2 \|\phi\|^2_{\dot H^s} \le  \frac{1}{2}   \|\phi\|^2_{\dot H^s}.
$$
Hence, we have $M^{(n)}_{c, \beta,m}  \ge c m^{2s}/2$.

Now let $\phi \in (H^s)^L $ with $\|\phi\|^2_2 =c$. We have
\begin{align}
I^{(n)}_{c, \beta,m}  (\phi_\lambda)    =    \frac{1}{2} \int (m^2+ |\lambda \xi|^2)^s |\widehat{\phi}(\xi)|^2 d\xi-  \lambda^{n-\beta}  \Upsilon_\beta (\phi).   \label{I10}
 \end{align}
We denote by $[s]$ the largest integer which is less than or equals to $s$, $\{s\}=s-[s]$. It suffices to consider the case that $s$ is not an integer.  Since
$$
(a+b)^s =(a+b)^{[s]} (a+b)^{\{s\}}  \leq   \sum^{[s]}_{j=0}\left(\begin{matrix}
j \\
 [s]
\end{matrix}\right)    a^j b^{( s -j)} +  \sum^{[s]}_{j=0} \left(\begin{matrix}
j \\
 [s]
\end{matrix}\right) a^{j+ \{s\}} b^{( [s] -j)},
$$
we have
\begin{align}
( m^2+ |\lambda \xi|^2)^s  & \leq  m^2+  \sum^{[s]}_{j=0} \left(\begin{matrix}
j \\
 [s]
\end{matrix}\right) m^{2j} (|\lambda \xi|^2)^{( s -j)} +  \sum^{[s]-1}_{j=0} \left(\begin{matrix}
j \\
 [s]
\end{matrix}\right) m^{2(j+ \{s\})} (|\lambda \xi|^2)^{( [s] -j)} \nonumber\\
& :=   m^2+  \lambda^{2\{s\}} P(\lambda,m, |\xi|). \label{I11}
 \end{align}
Noticing that for $\lambda \le 1$, we have $P(\lambda,m, |\xi|) \lesssim  1+|\xi|^2s$, which implied that
$$
\int \lambda^{2\{s\}} P(\lambda,m, |\xi|)|\widehat{\phi}(\xi)|^2 d\xi \to 0, \ \ \lambda \to 0.
$$
Hence, we have
\begin{align}
c m^{2s}/2 \leq I^{(n)}_{c, \beta,m}  (\phi_\lambda)   \leq c m^{2s}/2 + O(\lambda^{2\{s\}}),   \label{I12}
 \end{align}
which yields $M^{(n)}_{c, \beta,m} = c m^{2s}/2$. $\hfill\Box$

\begin{lem}
Let $s=(n-\beta)/2 > 1$. If $C^*c \le 1/2$, then $M^{(n)}_{c, \beta,m} $ is not achieved.
\end{lem}
{\bf Proof.}  Suppose on the contrary that there exists $u>0$ satisfying
$$
\frac{1}{2}c m^{2s} = I^{(n)}_{c, \beta,m}  (u) \ge   \frac{1}{2} \int (m^2+ | \xi|^2)^s |\widehat{u}(\xi)|^2 d\xi-   \frac{1}{2}  \|\phi\|^2_{\dot H^s}.
$$
 By the mean value theorem, there exits $\theta (t) \in (0,t)$ such that
$$
 f(t):= (m^2+t)^s -t^s -m^{2s} = st \left( (m^2+ \theta(t))^{s-1} - \theta (t)^{s-1}\right)>0
$$
for any $t>0$.  It follows that
 \begin{align}
 \frac{1}{2} cm^{2s} \geq \frac{1}{2} \int (m^2+ | \xi|^2)^s |\widehat{u}(\xi)|^2 d\xi-   \frac{1}{2}  \|\phi\|^2_{\dot H^s} =  \frac{1}{2} cm^{2s} + \int f(|\xi|^2) |\widehat{u}(\xi)|^2 d\xi.
 \end{align}
Noticing that $f(|\xi|^2)$ is a continuous functions of $\xi \in \mathbb{R}^n$ and $f(|\xi|^2)>0$ if $\xi \neq 0$, we immediately have $\int f(|\xi|^2) |\widehat{u}(\xi)|^2 d\xi>0$. A contraction. $\hfill \Box$

\medskip

Up to now, we have shown that for any $s>1$,  $I^{(n)}_{c, \beta,m}  (\cdot)$ has no minimizer.  In the following we consider the case $0<s<1$.

\begin{lem} \label{Ms<1}
Let $s=(n-\beta)/2 < 1$. If $C^*c < 1/2$, then $M^{(n)}_{c, \beta,m} \in (0,\  cm^{2s}/2).$
\end{lem}
{\bf Proof.} Let us denote $u_R := \mathscr{F}^{-1} \chi_{|\xi| \le R} \mathscr{F}u$.  Let $\phi \in (H^s)^L $ with $\|\phi\|^2_2= c$ satisfy
$$
\Upsilon_\beta (\phi)= C^*  \|\phi\|^2_2 \|\phi\|^2_{\dot H^s}.
$$
Then we have some $R>0$ satisfying
$$
\Upsilon_\beta (\phi_R ) \geq \frac{1}{2} C^* \|\phi\|^2_2  \|\phi\|^2_{\dot H^s} \geq \frac{1}{2} C^* \|\phi_R\|^2_2  \|\phi_R\|^2_{\dot H^s}.
$$
Taking $v= \sqrt{c} \phi_R/ \|\phi_R\|_2$, we see that $\|v\|^2_2=c$     and
$$
\Upsilon_\beta (v ) \geq \frac{1}{2} C^* \|v\|^2_2  \|v\|^2_{\dot H^s} .
$$
Moreover, the above inequality is invariant under the scaling $v \mapsto v_\lambda$, i.e.,
$$
\Upsilon_\beta (v_\lambda ) \geq \frac{1}{2} C^* \|v_\lambda\|^2_2  \|v_\lambda\|^2_{\dot H^s} = \frac{a}{2}   \|v_\lambda\|^2_{\dot H^s}, \ \ a= C^*c.
$$
Moreover, we have
\begin{align}
I^{(n)}_{c, \beta,m}  (v_\lambda)   \leq  \frac{1}{2} cm^{2s}  + \frac{1}{2} \int_{|\xi|\le R} \left( (m^2+ |\lambda \xi|^2)^s - a |\lambda\xi|^{2s} -m^{2s} \right)|\widehat{v}(\xi)|^2 d\xi.   \label{I18}
 \end{align}
Using the mean value theorem, for any $t>0$, we have some $\theta (t) \in (0,t)$ verifying
$$
f(t): = (m^2+t)^s - a t^s -m^{2s} = st \left( (m^2+ \theta(t))^{s-1} - a \theta (t)^{s-1}\right).
$$
Noticing that $s<1$, it follows that for $0<t\ll 1$, one has that
$$
 (m^2+ \theta(t))^{s-1} - a \theta (t)^{s-1} <0.
$$
Hence, taking $\lambda>0$ such that $\lambda R\ll 1$, we obtain that
$$
(m^2+ |\lambda \xi|^2)^s - a |\lambda\xi|^{2s} -m^{2s} <0, \ \ \forall \ 0<|\xi| \le R.
$$
Since $\xi \mapsto f(|\xi|^2)$ is continuous and $v\neq 0$, we immediately have $I^{(n)}_{c, \beta,m}  (v_\lambda) < cm^{2s}/2$.  Due to $C^*c <1/2$, we easily see that
$$
I^{(n)}_{c, \beta,m}  (\phi) > (1/2- C^*c)\|\phi\|^2_{\dot H^s} >0.
$$
The result follows. $\hfill \Box$

 \begin{lem} \label{Mc0}
Let $s=(n-\beta)/2 < 1$. If $C^*c = 1/2$, then $M^{(n)}_{c, \beta,m} =0.$
\end{lem}
{\bf Proof.}    Clearly, we have $M^{(n)}_{c, \beta,m}  \ge 0.$ Let us recall that for any minimizer $\phi$ of the functional $I^{(n)}_{c, \beta} (\cdot)$, we have for any $\varepsilon>0$,
\begin{align}
I^{(n)}_{c, \beta,m}  (\phi_\lambda)
&  =    \frac{1}{2} \int ((m^2+ |\lambda \xi|^2)^s -|\lambda\xi|^{2s}) |\widehat{\phi}(\xi)|^2 d\xi \nonumber\\
&  =    \frac{1}{2} \left(\int_{|\lambda \xi| > m/\varepsilon}+ \int_{|\lambda \xi| \le m/\varepsilon} \right) ((m^2+ |\lambda \xi|^2)^s -|\lambda\xi|^{2s}) |\widehat{\phi}(\xi)|^2 d\xi \nonumber\\
&:= I+II.
\label{I19}
 \end{align}
We estimate $I$. We may assume that $m/\varepsilon \gg 1$. Recall that
$$
(m^2+ |\lambda \xi|^2)^s -|\lambda\xi|^{2s} = |\lambda\xi|^{2s} \left( \left(1+ \frac{m^2}{|\lambda\xi|^{2 }}\right)^s -1\right) < s m^2  \frac{|\lambda\xi|^{2s } }{|\lambda\xi|^{2 }} \le  s\, m^{2s} \varepsilon^{2(1-s)}.
$$
It follows that
$$
I  \le c  s\, m^{2s} \varepsilon^{2(1-s)}.
$$
On the other hand, due to $(a+b)^s \le a^s +b^s$ and $\phi\in L^2$,
$$
II \le  \frac{1}{2} m^{2s} \int_{|\xi| \le m/\lambda \, \varepsilon} |\widehat{\phi}(\xi)|^2 d\xi \to 0, \ \ \lambda \to \infty.
$$
Hence, $I^{(n)}_{c, \beta,m}  (\phi_\lambda) \to 0$  as $\lambda\to \infty$. $\hfill\Box$

 \begin{lem} \label{1-critical}
Let $s=(n-\beta)/2 =1 $.   Then $M^{(n)}_{c, \beta,m} $ is achieved if and only if $C^*c=1/2$.
\end{lem}
{\bf Proof.} Noticing that for $s=1$
\begin{align}
I^{(n)}_{c, \beta,m}  (\phi )
&  =    \frac{1}{2} m^{2} c + \frac{1}{2} \int  | \xi|^2  |\widehat{\phi}(\xi)|^2 d\xi  - \Upsilon_\beta (\phi)= \frac{1}{2} m^{2} c + I^{(n)}_{c, \beta}  (\phi ),
\label{I29}
 \end{align}
we can obtain the result, as desired. $\hfill \Box$

By Lemma \ref{1-critical} and Theorem \ref{exist-minimizer2}, we can prove Theorem \ref{exist-minimizer4} in the case $s=1$.

\medskip

\noindent {\bf Proof of Theorem \ref{exist-minimizer4}.} In view of the discussions above, it suffices to consider the case $0<s<1$. Now let $u_k$ be a minimizing sequence. By  Lemma \ref{Ms<1}, we see that $u_k$ is bounded in $(H^s)^L $.  Following the proof as in Theorem \ref{exist-minimizer2}, we can assume that $u_k$ is radial and radially decreasing. We have
$$
M^{(n)}_{c, \beta,m} \le I^{(n)}_{c, \beta,m}  (u_k ) \to  M^{(n)}_{c, \beta,m}.
$$
Now we claim that $\inf \{\Upsilon_\beta (u_k) : \ k\ge 0\} \ge c_0$ for some $c_0>0$. If not, then we have $\Upsilon_\beta (u_k) \to 0$  up to a subsequence.  By Lemma \ref{Ms<1},
$$
\frac{1}{2} m^{2s} c \le \lim_{k\to \infty} \frac{1}{2} \|(m^2+|\xi|^2)^{s/2} \widehat{u}_k\|^2_2 = \lim_{k\to \infty}  I^{(n)}_{c, \beta,m}  (u_k ) =  M^{(n)}_{c, \beta,m} < \frac{1}{2} m^{2s} c.
$$
This is a contradiction.

Now we can repeat the same procedure as in the proof of Theorem \ref{exist-minimizer2} to show that $u_k \to u \ge 0$ and $u\neq 0$, with a minimizer $u$, as desired.

Finally, we show the necessity of $C^*c<1/2$.  If not,
then $C^*c\geq 1/2$. If $C^*c>1/2$, be Lemma \ref{Mcinfty} we have $M^{(n)}_{c, \beta,m}=-\infty$. If $C^*c=1/2$, in view of Lemma \ref{Mc0} we have   $M^{(n)}_{c, \beta,m} =0 $. If  $u\neq 0$ is a minimizer, then $I^{(n)}_{c, \beta,m}  (u) = 0$. On the other hand,  from the definition of $I^{(n)}_{c, \beta,m}  (\cdot)$ we have  $I^{(n)}_{c, \beta,m}  (u) > 0$. A contradiction.
 $\hfill \Box$\\

 \begin{appendix}
 \section{Proof of Theorem \ref{Prop4.3.1}}

The proof of Theorem \ref{Prop4.3.1} is essentially known and we now sketch its proof  by  following  \cite{WHHG}, Section 2.4 (see also \cite{KeKo09} in 3D).
 \begin{prop}\label{Propsem}
 Let $H(t)=e^{t\Delta}$, $\mathscr{A} f =\int^t_0 H(t-s) f(s) ds$. We have
\begin{align}
&  \|H(t)u_0\|_{L^{n+2}_{x, t\in [0,T]}} \lesssim
\|u_0\|_{n},
\label{eq4.3.6}\\
&  \| H(t)u_0\|_{L^{\infty}(0,T; \ L^{n})} \lesssim \|u_0\|_{n},
\label{eq4.3.7}\\
 &  \|\nabla \mathscr{A} f\|_{L^{n+2}_{x, t\in [0,T]}}
\lesssim \|f\|_{L^{(n+2)/2}_{x, t\in [0,T]} },
\label{eq4.3.8}\\
&  \| \nabla \mathscr{A} f\|_{ L^\infty (0,T; \ L^n) } \lesssim
\|f\|_{L^{(2+n)/2}_{x, t\in [0,T]}}. \label{eq4.3.9}
\end{align}
\end{prop}
Put
\begin{align}
& \mathfrak{D} =\left\{u: \|u\|_{L^{2+n}_{x, t\in [0,T]}} \le
\delta, \   \|
u\|_{L^\infty ([0,T]; L^n) } \le 2 C \|u_0\|_n \right\}, \label{eq4.3.15}\\
& d(u,v) =  \|u-v\|_{L^{2+n}_{x, t\in [0,T]}}. \label{eq4.3.16}
\end{align}
We consider the mapping:
\begin{align}
\mathfrak{M}:  u(t) \to H(t) u_0 + \mathscr{A} \mathbb{P} \, {\rm
div}\, (u \otimes u),  \label{eq4.3.17}
\end{align}
where
\begin{align}\label{eq4.0.4}
\mathbb{P} = I + (-\Delta)^{-1} \nabla  {\rm div} .
\end{align}
By Proposition \ref{Propsem}, we have
\begin{align}
 \|\mathfrak{M}  u\|_{L^{2+n}_{x, t\in [0,T]}}  & \lesssim  \|H(t) u_0
\|_{L^{2+n}_{x, t\in [0,T]} }+ \| u \otimes  u \|_{L^{(2+n)/2}_{x,
t\in [0,T]}} \nonumber\\
 & \lesssim  \|H(t) u_0
\|_{L^{2+n}_{x, t\in [0,T]} }+ \delta^2, \label{eq4.3.18} \\
 \|\mathfrak{M}  u\|_{L^{\infty}(0,T; \  L^n)}  & \lesssim  \|u_0
\|_{n}+ \| u \otimes  u \|_{L^{(2+n)/2}_{x,
t\in [0,T]}} \nonumber\\
 & \lesssim  \|u_0
\|_{n}+  \delta^2. \label{eq4.3.19}
\end{align}
If $C\delta \le 1/4$,  we can show that $\mathfrak{M}$
is a contraction mapping from $ \mathfrak{D}$ into itself.  So, there exists a $u$ satisfying
\begin{align}
 u(t) = H(t) u_0 + \mathscr{A} \mathbb{P} \nabla \cdot (u \otimes u).
\label{eq4.3.20}
\end{align}
By a standard argument, we see that $u$ is unique in $L^{2+n}(0,T; \ L^{2+n} )$.
Moreover, one can extend the solution step by step and find a maximal
$T_m$ such that  $u \in C ([0,T_m); \  L^n )  \cap
L^{2+n}_{\rm loc}(0,T_m; \  L^{2+n} )$.  In the following we show that
$$
\|u\|_{L^{2+n} (0,T_m; \  L^{2+n} )}=\infty.
$$
Assume for a contrary that $\|u\|_{L^{2+n} (0,T_m; \  L^{2+n} )}<\infty$. In view of the first inequality in \eqref{eq4.3.19} we see that
$$
\|u\|_{C ([0,T_m); \  L^n ) \ \cap \ L^{2+n} (0,T_m; \  L^{2+n} )}<\infty.
$$
Using the same idea as in \cite{CaWe89} for the nonlinear Schr\"odinger equation,
we now extend the solution beyond $T_m$. We have for $0< T_m-T \ll 1$,
$$
  u(t)= H(t-T)u(T) +  \int^t_{T} H(t-\tau )\mathbb{P} \nabla \cdot (u \otimes u)d\tau.
$$
It follows that
\begin{align}
 \|H     (t-T) u(T)  \|_{L^{n+2}_{x,t\in  (T, T_m)}}   
    \le  \|u\|_{L^{n+2}_{x,t\in (T, T_m)}}+
  \left\|   u  \right\|^2_{L^{n+2}_{x,t\in (T, T_m)}} \to 0, \ \ T\to T_m.
\end{align}
Replacing $[0,T]$ by $[T, T_m]$ and $\|u_0\|_n$ by $\|u(T)\|_n$ in the definition of $(\mathfrak{D},d)$,  we can find that  the solution can be extended to $C([T,T_m],L^n)$ if $T$ is sufficiently close to $T_m$. It follows that the solution exists beyond $T_m$. A contradiction. \\

\end{appendix}

{\bf Acknowledgment.} Part of the work was carried out while the fourth named author was visiting the LMPT at Universit\'e Fran\c cois Rabelais Tours. He is grateful to LMPT for its hospitality.

\end{document}